\UseRawInputEncoding
\documentclass[12pt]{article}

\usepackage{float}
\usepackage{latexsym}
\usepackage{mathrsfs}
\usepackage{amsfonts,amsmath,amssymb}
\usepackage{subfigure}
\usepackage[ruled]{algorithm2e}
\usepackage{algpseudocode}
\usepackage{indentfirst}
\usepackage{tabularx}
\usepackage{booktabs}
\usepackage{longtable}
\usepackage{multirow}
\usepackage{caption}
\usepackage{subeqnarray}
\usepackage{verbatim}
\usepackage{tikz}

\newcommand{\smallcircle}{
  \tikz \draw[fill=black] (0,0) circle (1.5pt);
}

\newtheorem{theorem}{Theorem}[section]
\newtheorem{lemma}[theorem]{Lemma}

\newtheorem{definition}[theorem]{Definition}

\numberwithin{equation}{section}
\newenvironment{proof}[1][Proof]{\textbf{\varphi1.} }
{\ \rule{0.75em}{0.75em}\smallskip}

\usepackage[pdftex,unicode,colorlinks,
linkcolor=blue]{hyperref}

\def\bar{\overline}

\def\lista\varphi1
{{ \itemindent 0.0cm \labelsep .2cm \leftmargin 0.8cm \rightmargin
		0.0cm \labelwidth 0.6cm \topsep 0.0mm
		\parsep 0.0mm
		\itemsep 0.0mm
		\begin{list}{}
			{ \setlength{\leftmargin}{.8cm} \setlength{\rightmargin}{0.0cm}
				\setlength{\parsep}{0.0mm} \setlength{\topsep}{.0mm}
				\setlength{\parskip}{.0cm} \setlength{\itemsep}{.0cm} }
			{\varphi1}\end{list}} }

\begin{document}
	
\title{\Large \bf
An inexact golden ratio primal-dual algorithm with  linesearch step for a saddle point problem}

\author{
Changjie Fang \,\footnote{\, College of Science, Chongqing University of Posts and Telecommunications,	Chongqing 400065, China.
			E-mail: {\tt fangcj@cqupt.edu.cn.}},
		\ \
		Jinxiu Liu \,\footnote{\, College of Science, Chongqing University of Posts and Telecommunications, Chongqing 400065, China.
			E-mail: {\tt liujxtm@163.com.}},
\ \
	Jingtao Qiu \,\footnote{\, College of Science, Chongqing University of Posts and Telecommunications, Chongqing 400065, China. E-mail:{\tt qiujingtao007@126.com.}}
\ \ \mbox{and} \\
Shenglan Chen \,\footnote{\, College of Science, Chongqing University of Posts and Telecommunications, Chongqing 400065, China. E-mail:{\tt chensl@cqupt.edu.cn.}}
}

\date{}
\maketitle

\noindent {\bf Abstract}
In this paper, we propose an inexact golden ratio primal-dual algorithm with linesearch step(IP-GRPDAL) for  solving the saddle point problems,  where two subproblems can be approximately solved by applying the notations of inexact extended proximal operators with  matrix norm. Our proposed IP-GRPDAL method allows for larger stepsizes by replacing the extrapolation step with a convex combination step. Each iteration of the linesearch requires to update only the dual variable, and hence it is quite cheap. In addition, we prove convergence of the proposed algorithm and show an $O(1/N)$ ergodic convergence rate for our algorithm, where $N$ represents the number of iterations. When one of the component functions is strongly convex, the accelerated $O(1/N^2)$ convergence rate results are established by choosing adaptively some algorithmic parameters. Furthermore, when  both component functions are strongly convex, the linear convergence rate results are achieved. Numerical simulation results on the sparse recovery and image deblurring problems illustrate the feasibility and efficiency of our inexact algorithms.
\vspace{5mm}	

\vspace{5mm}

\section{Introduction}
\noindent
 Let $X:=\mathbb{R}^{n}$ and $Y:=\mathbb{R}^{m}$ be two finite-dimensional Euclidean spaces equipped with a standard inner product $\langle\cdot,\cdot\rangle$ and a norm $\|\cdot\|=\sqrt{\langle\cdot,\cdot\rangle}$.
	Let $f: X \to (-\infty,+\infty]$ and $g,h: Y \to (-\infty,+\infty]$ be proper lower semicontinuous (l.s.c) convex functions,  $A:X \to Y$ be a bounded linear mapping. Denote  the Legendre-Fenchel conjugate of  $h$ and the adjoint of  $A$ by $h^*$ and $A^*$, respectively. Now we consider the primal problem
	\begin{align} \label{1.3}
		\mathop{\min}\limits_{x \in X } f(x) + h(Ax)
	\end{align}
	together with its dual problem
	\begin{align}\label{1.3a}
	 \mathop{\max}\limits_{y \in Y } -f^*(-A^*y) - h^*(y).
	\end{align}
	If a primal-dual solution pair $(\overline{x},\overline{y})$ of \eqref{1.3} and \eqref{1.3a} exists, i.e.,
\begin{align*}
 0\in\partial f(\overline{x})+A^*\overline{y},\,\, 0\in\partial h(A\overline{x})-\overline{y},
\end{align*}
then the problem \eqref{1.3} is equivalent to the following saddle-point formulation:
\begin{align}\label{1.3b}
  \min_{x\in X}\max_{y\in Y}f(x)+\langle Ax,y\rangle-h^*(y).
\end{align}

  It is well known that many application problems can be formulated as the saddle point problem (\ref{1.3b}) such as image restoration, magnetic resonance imaging and computer vision; see, for example, \cite{PT2009,RL1992,ZC2008}. Two of the most popular approaches are the primal-dual algorithm (PDA) \cite{CP2011,HB2012}, alternating direction method of multipliers (ADMM) method \cite{BS2010,HBS2012}, and their accelerated and generalized variants \cite{LY2018,YM2018,YM}. To solve model (\ref{1.3b}), the following first-order primal-dual algorithm (PDA) \cite{CP2011} has attracted much attention:
	\begin{equation} \label{1.5}
		\begin{cases}
			{x}^{k+1} = \textrm{Prox}_{\tau f}(x^k-\tau A^*y^k),\\
			\hat{x}^{k+1} = x^{k+1}+\theta({x^{k+1}-x^k}),\\
			{y}^{k+1} = \textrm{Prox}_{\sigma h^*}(y^k+\sigma Ax^{k+1}),
		\end{cases}
	\end{equation}
where $\tau>0$ and $\sigma>0$ are regularization parameters and $\theta \in (0,1)$ is an extrapolation parameter, and for $\theta =1$, the convergence of PDA was proved with the requirement on step sizes $\tau\sigma \|A\|^2<1$. Generally, with fixed $\tau$ and $\sigma$, a flexible extrapolation parameter $\theta$ is of benefit to a potential acceleration of the algorithm  (\ref{1.5}), which motivates researchers to enlarge the range of $\theta$, e.g., see \cite{HJ2016,KH2020}. Indeed, the scheme (\ref{1.5}) reduces to the classic Arrow-Hurwicz method \cite{Arrow1958} when $\theta=0$. However, the convergence of the Arrow-Hurwicz method can only be guaranteed under the restrictive assumption that $\tau$ and $\sigma$ are sufficiently small. In order to overcome this difficulty, Chang et al. in \cite{CY} proposed a golden ratio primal-dual algorithm (GRPDA) for solving \eqref{1.3b}  based on a seminal convex combination technique introduced by Malitsky in \cite{YM}, which,  started at $(x^0,y^0)\in\mathbb{R}^{n}\times\mathbb{R}^{m} $, iterates as
	\begin{equation}
		\begin{cases}
			z^{k+1}=\frac{\phi-1}{\phi}x^k+\frac{1}{\phi}z^k,\\
			{x}^{k+1} = \textrm{Prox}_{\tau f}(z^{k+1}-\tau A^*y^k),\\
			{y}^{k+1} = \textrm{Prox}_{\sigma h^*}(y^k+\sigma Ax^{k+1}),
		\end{cases}
	\end{equation}
where $\phi\in (1,\frac{1+\sqrt{5}}{2}]$ determines the convex combination coefficients. In \cite{CY}, the iterative convergence and ergodic convergence rate results are established under the condition $\tau\sigma \|A\|^2 < \phi$. Since $\phi>1$, this stepsize condition is much relaxed than that of the PDA method \eqref{1.5}; see also \cite{CY2022}. Further, Chang et al. incorporated linesearch strategy into the GRPDA method, in which the next iteration $y^{k+1}$ is implicitly computed, since the stepsize parameter $\tau_{k+1}$ is determined by the inequality including $y^{k+1}$; see Step 2  of  Algorithm 3.1 in \cite{C2022}.

As primal-dual algorithms, however, when the proximal operators of  $f$ and $h^*$  are not easy to compute, they did not perform ideally well in terms of computing time and efficiency, e.g., see the examples in \cite{CC2010,EB2016} and the numerical experiments in \cite{HH2014}. And in many practical applications, one often encounters the case that at least one of the proximal operators does not possess a closed-form solution and their evaluation involves inner iterative algorithms.  In this situation, some researchers are dedicated to approximately solving the subproblems instead of finding their accurate solutions, for example, \cite{HZ2024,LY2018,NM2011,X2018,WS2024}. An absolute error criterion was adopted in \cite{EJ1992}, where the subproblem errors are controlled by a summable sequence of error tolerances.  Jiang et al. \cite{FC2021,FW2021} studied  two inexact primal-dual algorithms with absolute and relative error criteria respectively, where for the inexact primal-dual method with a relative error criterion, only the $O(1/N)$ convergence rate was established.  In \cite{RJ2020}, Rasch and Chambolle proposed the inexact first-order primal-dual methods by applying the concepts of inexact proxima where all the controlled errors were required to be summable. Further, Fang et al. \cite{FH2022} proposed an  inexact primal-dual method with correction step by introducing the notations of inexact extended proximal operators with  matrix norm, where the $O(1/N)$ ergodic convergence rate was achieved. In \cite{FH2022}, the accelerated versions of the proposed method under the assumptions that $f$ or $h^*$  is strongly convex have not been considered.

  In this paper, we are concerned with the following saddle point problem:
	\begin{align}\label{1.1}
		\mathop{\min}\limits_{x \in X }\mathop{\max}\limits_{y \in Y}L(x,y)=f(x)+{\langle{Ax,y}\rangle}-g(y)
	\end{align}
	 Recall that $(\bar{x},\bar{y})$ is called a saddle point of (\ref{1.1}) if it satisfies the inequalities
	\begin{align}
		L(\bar{x},y) \leqslant L(\bar{x},\bar{y}) \leqslant L(x,\bar{y}),  \forall x \in  X,  \forall y \in Y
	\end{align}
Hence, Problem \eqref{1.3b} is a special case of Problem \eqref{1.1}.

	Motivated by the research works \cite{C2022,FH2022}, in this paper, we propose an inexact golden ratio primal-dual algorithm with  linesearch step for solving problem \eqref{1.1} by applying the type-2 approximation of the extended proximal point introduced in \cite{FH2022}. The main contributions of this paper are summarized as follows.\par
\textbullet \quad For the case  the proximal operators of  $f$ and $g$  are not easy to compute,  we propose an inexact IP-GRPDAL method with extended proximal terms containing symmetric positive definite matrix. Both subproblems in our  method can be approximately solved under the type-2 approximation criterias. Global convergence and $O(1/N)$ ergodic convergence rate results  are established under the condition $\tau\sigma \|A\|_T^2 < \phi$, where $N$ denotes the iteration counter. Furthermore, we establish the convergence rates in case the error tolerances $\{\delta_k\}$ and $\{\varepsilon_k\}$ are required to decrease like $O(1/k^{2\alpha+1})$ for some $\alpha > 0$. \par
\textbullet \quad Our method updates the dual variable by adopting linesearch step to allow adaptive and potentially much larger stepsizes which effectively reduces the computational effort of the algorithm iteration. In addition, the next iteration $y^{k+1}$ is explicitly computed compared with that in \cite{C2022}; see \eqref{3.3} in Algorithm \ref{A1}. \par
\textbullet \quad We propose the accelerated versions of IP-GRPDAL method, which were not provided in \cite{FH2022}. When one of the  underlying functions is strongly convex, $O(1/N^2)$ convergence rate results are established by adaptively choosing some algorithmic parameters, for example, $\beta$ is replaced by $\beta_k$; see \eqref{3.32} of Algorithm \ref{A2}. In addition, the linear convergence rate results can be established when both $f$ and $g$ are strongly convex. \par
 \textbullet \quad We perform numerical experiments for the sparse recovery and image deblurring problems, demonstrating that our method outperforms some existing methods \cite{C2022,CP2011,FH2022,YM2018}.\par
	The rest of this paper is organized as follows. In Section \ref{sec:pre}, we introduce the concepts of inexact extended proximal terms and present some auxiliary material. In Section \ref{sec:alg}, the main  algorithm and  its accelerated versions are presented. At the same time, we also prove the convergence of our algorithms and analyze their convergence rates. Numerical experiment results are reported in Section \ref{sec:exp}. Some conclusions are presented in Section \ref{sec:clu}.

\section{Preliminaries}\label{sec:pre}
	For given $x_1 \in X$ and $y_1 \in Y$, we define $P_{x_1,y_1}(x) := f(x)-f(x_1)+{\langle{x-x_1,A^*y_1}\rangle}$ and $D_{x_1,y_1}(y) := g(y)-g(y_1)+{\langle{y-y_1,-Ax_1}\rangle}$ for any $x \in  X$ and $y \in Y$. Let $(\bar{x},\bar{y}) \in X \times Y$ be a generic saddle point. When there is no confusion, we will omit the subscript in $P$ and $D$,
	\begin{equation} \label{2.1}
		\begin{cases}
			P(x) := P_{\bar{x},\bar{y}}(x) = f(x) - f(\bar{x}) + \langle x-\bar{x}, A^*\bar{y} \rangle, \forall x \in X,\\
			D(y) := D_{\bar{x},\bar{y}}(y) = g(y) - g(\bar{y}) + \langle y-\bar{y}, -A\bar{x} \rangle, \forall y \in Y.
		\end{cases}
	\end{equation}
	By subgradient inequality, it is clear that $P(x) \geqslant 0$ and $D(y) \geqslant 0$. Note that the functions $P(x)$ and $D(y)$ are convex in $x$ and $y$, respectively. The primal-dual gap is defined as $G(x,y) := L(x,\bar{y})-L(\bar{x},y)$ for $(x,y) \in X \times Y$. It is easy to verify that
	\begin{equation}
		G(x, y) := G_{\bar{x},\bar{y}}(x, y) = P(x) + D(y) \geqslant 0, \forall (x,y) \in X \times Y.
	\end{equation}
	The system (\ref{2.1}) can be reformulated as
	$$ 0 \in \partial f(\bar{x})+A^*\bar{y},\qquad 0 \in \partial g(\bar{y})-A\bar{x}. $$
	 Suppose that $h$ is a convex function in $\mathbb{R}^{n}$, and that $D \in \mathbb{R}^{n \times n}$ is a symmetric positive definite matrix. For any $D \succ 0$ and given $y\in\mathbb{R}^{n}$, denote
	\begin{equation}
		J_y(x) := h(x) + \frac{1}{2\tau} {\|x-y\|}^2_D, \forall x \in \mathbb{R}^{n},
	\end{equation}
	and define the proximal operator of $h$ as
	\begin{equation}
		\textnormal{Prox}_{\tau h}^D(y) = \mathop{\arg\min}\limits_{x \in X } \{h(x) + \frac{1}{2\tau} {\|x-y\|}^2_D\},
	\end{equation}
	where ${\|x\|}_D^2 = \langle{x,Dx}\rangle$ and $D^{-1}$ denotes the inverse of $D$, the first-order optimality condition for the proximum gives different characterizations of the proximal operator:
	$$ \bar{z} :=  \textrm{Prox}_{\tau h}^D(y) \Longleftrightarrow 0 \in \partial J_y(\bar{z}) \Longleftrightarrow \frac{y-\bar{z}}{\tau} \in \partial h(\bar{z}).$$
	Below, we recall  definitions of three different types of inexact  extended proximal operators with  matrix norm, which can be found in \cite{FH2022}.
	
	\begin{definition} Let $\varepsilon \geqslant 0$. $z \in X$ is said to be a type-0 approximation of the extended proximal point $Prox_{\tau h}^D(y)$ with precision $\varepsilon$ if
		\begin{equation}
			z \thickapprox_0^{\varepsilon} \textnormal{Prox}_{\tau h}^D(y) \Longleftrightarrow \|z-\bar{z}\|_D \leqslant \sqrt{2\tau\varepsilon}
	\end{equation}
\end{definition}

\begin{definition} Let $\varepsilon \geqslant 0$. $z \in X$ is said to be a type-1 approximation of the extended proximal point $Prox_{\tau h}^D(y)$ with precision $\varepsilon$ if
	\begin{equation}
		z \thickapprox_1^{\varepsilon}  \textnormal{Prox}_{\tau h}^D(y) \Longleftrightarrow 0 \in \partial_\varepsilon J_y(z),
	\end{equation}
	where $ \partial_\varepsilon J_y(z) = \{p \in X | J_y(x) \geqslant  J_y(z) + \langle{p,x-z}\rangle - \varepsilon, \forall x \in X \}.$
\end{definition}

\begin{definition} \label{D2.3} Let $\varepsilon \geqslant 0$. $z \in X$ is said to be a type-2 approximation of the extended proximal point $Prox_{\tau h}^D(y)$ with precision $\varepsilon$ if
	\begin{equation}
		z \thickapprox_2^{\varepsilon}\textnormal{Prox}_{\tau h}^D(y) \Longleftrightarrow \frac{1}{\tau}D(y-z) \in \partial_\varepsilon h(z),
	\end{equation}
	where $ \partial_\varepsilon h(z) = \{p \in X | h(x) \geqslant  h(z) + \langle{p,x-z}\rangle - \varepsilon, \forall x \in X \}.$
\end{definition}

We give two simple but useful lemmas in the following.
\begin{lemma}
	For any $x,y,z$ and a symmetric positive definite matrix $D$, we have the identity
	\begin{equation} \label{2.8}
		\langle{D(x-y),x-z}\rangle = \frac{1}{2}[\|x-y\|_D^2 + \|x-z\|_D^2 - \|y-z\|_D^2].
	\end{equation}
	For $\alpha \in\mathbb{R}$, there holds
	\begin{equation} \label{2.9}
		\|\alpha x + (1-\alpha)y\|_D^2 = \alpha \|x\|_D^2 + (1-\alpha) \|y\|_D^2 - \alpha(1-\alpha)\|x-y\|_D^2.
	\end{equation}
\end{lemma}

\begin{lemma}\cite{Shores2007}
	Assuming that $\lambda$ and $\Lambda$ are the minimum and maximum eigenvalues of the symmetric positive definite matrix $D$, respectively, we have
	\begin{equation} \label{2.10}
		\sqrt{\lambda} \|x\| \leqslant \|x\|_D \leqslant \sqrt{\Lambda} \|x\|.
	\end{equation}
\end{lemma}

\section{Algorithm and convergence properties}\label{sec:alg}
In this section, we propose an inexact GRPDA algorithm and then show the convergence of the proposed method. If $f$ is further assumed to be strongly convex, we can modify our method to accelerate the convergence rate. Moreover, if both $f$ and $g$ are strongly convex, a linear convergence rate can be achieved.
\subsection{Convex case}
	\begin{algorithm}[H] \label{A1}
		 1:\quad Let $\varphi=\frac{\sqrt{5}+1}{2}$ be the golden ratio, that is $\varphi^2=1+\varphi$. Choose $ x^0 = z^0 \in 	\mathbb{R}^{n}, y^0 \in \mathbb{R}^{m}, \phi \in (1,\varphi), \eta \in (0,1), \mu \in (0,1)$, $ \tau_0 >0$ and $\beta > 0$. $S \in \mathbb{R}^{n \times n}$ and $T \in \mathbb{R}^{m \times m}$ are given symmetric positive definite matrix. Set $\psi = \frac{1+\phi}{\phi^2}$ and $ k=0.$\\
		2:\quad  Compute
		\begin{equation} \label{3.1}
			{z}^{k+1} = \frac{\phi-1}{\phi} x^k + \frac{1}{\phi} z^k,
		\end{equation}
		\begin{equation} \label{3.2}
			{x}^{k+1} \approx_2^{\delta_{k+1}} \mathop{\arg\min}\limits_{x \in X} \{ L(x,y^k) + \frac{1}{2\tau_k}\|x-z^{k+1}\|_S^2 \}.
		\end{equation}\\
		3: \quad  Choose any $ \tau_{k+1} \in [\tau_{k},\psi\tau_{k}]$ and run \\
		\quad \quad 3.a: Compute
		\begin{equation} \label{3.3}
			{y}^{k+1} \approx_2^{\varepsilon_{k+1}} \mathop{\arg\max}\limits_{y \in Y} \{ L(x^{k+1},y) - \frac{1}{2\beta \tau_{k+1}}\|y-y^k\|_T^2 \}.
		\end{equation}
		\quad \quad 3.b: Break linesearch if
		\begin{equation} \label{3.4}
			\sqrt{\beta\tau_{k+1}}\|A^*y^{k+1}-A^*y^k\|_T \leqslant \eta\sqrt{\frac{\phi}{\tau_k}} \|y^{k+1}-y^k\|_T.
	\end{equation}
	\quad \quad Otherwise, set $\tau_{k+1} := \tau_{k+1}\mu$ and go to 3.a.\\
   4:\quad Set $k+1 \leftarrow k+2$ and return to 2.
	\caption{An inexact GRPDA with linesearch(IP-GRPDAL)}
\end{algorithm}

Firstly, we summarize several useful lemmas  which will be used in the sequence..
\begin{lemma} \label{L3.1}
(i) The linesearch in Algorithm \ref{A1} always terminates.\\
(ii) There exists $\underline{\tau} := \frac{\eta\sqrt{\phi}}{L\sqrt{\beta\psi}} > 0 $ such that $\tau_{k+1} > \underline{\tau}$ for all $k \geqslant 0$. \\
(iii) For any integer $N>0$, we have $|\mathcal{A}_N| \geqslant \hat{c}N$ for some constant $\hat{c}>0$, where $\mathcal{A}_N = \{1 \leqslant k \leqslant N:\tau_{k} \geqslant \underline{\tau}\}$ and $|\mathcal{A}_N|$ is the cardinality of the set $\mathcal{A}_N$, which implies $\sum_{k=1}^{N} \tau_{k+1} \geqslant \underline{c}N$ with $\underline{c} = \hat{c}\underline{\tau}$.
\end{lemma}
\noindent Proof. (i) In each iteration of the linesearch $\tau_{k}$ is multiplied by factor $\mu \in (0,1)$. Since (\ref{3.4}) is fulfilled whenever $\tau_{k} \leqslant \frac{\eta\sqrt{\phi}}{L\sqrt{\beta\psi}}$ where $L=\|A^*\|_T$, the inner loop can not run indefinitely.

(ii)According to the recursive method, we assume that $\tau_0 > \frac{\eta\sqrt{\phi}}{L\sqrt{\beta\psi}}$. Then, our goal is to show that from $\tau_{k} > \frac{\eta\sqrt{\phi}}{L\sqrt{\beta\psi}} $ follows $\tau_{k+1} > \frac{\eta\sqrt{\phi}}{L\sqrt{\beta\psi}}$. Suppose that $\tau_{k+1} = \psi\tau_{k}\mu^i $ for some $i \in Z^+$. If $i=0$ then $ \tau_{k+1} > \tau_{k} > \frac{\eta\sqrt{\phi}}{L\sqrt{\beta\psi}} $. If $i>0$ then $\tau_{k+1}^{'} = \psi\tau_{k}\mu^{i-1}$ does not satisfy (\ref{3.4}). Thus, $\tau_{k+1}^{'} > \frac{\eta\sqrt{\phi}}{L\sqrt{\beta\psi}}$ and hence, $\tau_{k+1} > \frac{\eta\sqrt{\phi}}{L\sqrt{\beta\psi}} $.

(iii) The detailed proof takes similar approach with Lemma 3.1 (iii) of $\cite{C2022}$ and is thus omitted.

\begin{lemma} \label{L3.2}
Suppose that $\lambda_1, \lambda_2 > \eta $ are the minimum eigenvalues of $S$ and $T$, respectively. Let $\theta_{k+1} = \frac{\tau_{k+1}}{\tau_k}$ and $\{(z^{k+1}, x^{k+1}, y^{k+1}) : k \geqslant\ 0\}$ be the sequence generated by Algorithm \ref{A1}. Then, for any $ (\bar{x}, \bar{y}) \in X \times Y$, there holds
\begin{equation} \label{3.5}
	\begin{split}
		\tau_{k+1}&G(x^{k+1}, y^{k+1}) \leqslant  \langle S(x^{k+2} - z^{k+2}), \bar{x}-x^{k+2} \rangle +\frac{1}{\beta}\langle T(y^{k+1} - y^k), \bar{y}-y^{k+1} \rangle  \\
		&+\phi\theta_{k+1} \langle S(x^{k+1} - z^{k+2}), x^{k+2}-x^{k+1} \rangle + \tau_{k+1} \langle A^*(y^{k+1}-y^{k}), x^{k+1}-x^{k+2}\rangle\\
		&+\tau_{k+1}(\delta_{k+1}+\delta_{k+2}+\varepsilon_{k+1}).
	\end{split}
\end{equation}
\end{lemma}

\noindent Proof. By Definition 2.3, the optimal condition of (\ref{3.2}) yields
$$\frac{1}{\tau_k} S(z^{k+1}-x^{k+1}) - A^*y^k \in \partial_{\delta_{k+1}} f(x^{k+1}).$$
In view of the definition of $\varepsilon$-subdifferential, we have
\begin{equation} \label{3.6}
f(x)-f(x^{k+1})+ \frac{1}{\tau_{k}}\langle S{(x^{k+1}-z^{k+1}) + \tau_{k}A^*y^k, x-x^{k+1}} \rangle + \delta_{k+1} \geqslant 0, \forall x \in X.
\end{equation}
Setting $x=\bar{x}$ and $x=x^{k+2}$ in (\ref{3.6}), respectively, we obtain
\begin{equation} \label{3.7}
\tau_{k}(f(x^{k+1})-f(\bar{x})) \leqslant  \langle S(x^{k+1}-z^{k+1}) + \tau_{k}A^*y^k, \bar{x}-x^{k+1} \rangle + \tau_{k}\delta_{k+1}.
\end{equation}
\begin{equation} \label{3.8}
\tau_{k}(f(x^{k+1})-f(x^{k+2}) \leqslant \langle S(x^{k+1} - z^{k+1}) + \tau_{k}A^*y^{k}, x^{k+2}-x^{k+1} \rangle + \tau_{k}\delta_{k+1},
\end{equation}
In view of (\ref{3.7}), we have
\begin{equation} \label{3.9}
\tau_{k+1}(f(x^{k+2})-f(\bar{x})) \leqslant \langle S(x^{k+2} - z^{k+2}) + \tau_{k+1}A^*y^{k+1}, \bar{x}-x^{k+2} \rangle + \tau_{k+1}\delta_{k+2}.
\end{equation}
Multiplying (\ref{3.8}) by $\theta_{k+1}$ and using the fact $x^{k+1} -z^{k+1} = \phi(x^{k+1}-z^{k+2})$, we obtain
\begin{equation} \label{3.10}
\tau_{k+1}(f(x^{k+1})-f(x^{k+2})) \leqslant \langle \phi\theta_{k+1} S(x^{k+1} - z^{k+2}) + \tau_{k+1}A^*y^{k}, x^{k+2}-x^{k+1} \rangle + \tau_{k+1}\delta_{k+1}.
\end{equation}
Similarly, for $\bar{y} \in Y$, from (\ref{3.3}) we obtain
\begin{equation} \label{3.11}
\tau_{k+1}(g(y^{k+1})-g(\bar{y})) \leqslant \frac{1}{\beta} \langle T(y^{k+1} - y^k) - \beta\tau_{k+1}Ax^{k+1}, \bar{y}-y^{k+1} \rangle + \tau_{k+1}\varepsilon_{k+1}.
\end{equation}
Direct calculations show that a summation of (\ref{3.9}), (\ref{3.10}) and (\ref{3.11}) gives
$$ \tau_{k+1}(f(x^{k+1})-f(\bar{x})) + \tau_{k+1}(g(y^{k+1})-g(\bar{y})) + \tau_{k+1} \langle A^*\bar{y}, x^{k+1}-\bar{x}\rangle - \tau_{k+1} \langle A\bar{x}, y^{k+1}-\bar{y} \rangle $$
$$  \leqslant  \langle S(x^{k+2} - z^{k+2}), \bar{x}-x^{k+2} \rangle + \tau_{k+1}\langle A^*y^{k+1}, \bar{x} -x^{k+2} \rangle + \phi\theta_{k+1} \langle S(x^{k+1} - z^{k+2}), x^{k+2}-x^{k+1} \rangle $$
$$+ \tau_{k+1}\langle A^*y^k, x^{k+2}-x^{k+1} \rangle + \frac{1}{\beta}\langle T(y^{k+1} - y^k), \bar{y}-y^{k+1} \rangle  - \tau_{k+1} \langle Ax^{k+1}, \bar{y} - y^{k+1} \rangle $$
$$ + \tau_{k+1} \langle A^*\bar{y}, x^{k+1}-\bar{x}\rangle -  \tau_{k+1} \langle A\bar{x}, y^{k+1}-\bar{y}\rangle+ \tau_{k+1}(\delta_{k+1}+\delta_{k+2}+\varepsilon_{k+1}).$$
Applying the definitions of $P(x)$ and $D(y)$ in (\ref{2.1}) to  the above inequality, we have
$$ \tau_{k+1}P(x^{k+1}) + \tau_{k+1}D(y^{k+1}) \leqslant \langle  S(x^{k+2} - z^{k+2}), \bar{x}-x^{k+2} \rangle + \phi\theta_{k+1} \langle S(x^{k+1} - z^{k+2}), x^{k+2}-x^{k+1} \rangle $$
$$ \frac{1}{\beta}\langle T(y^{k+1} - y^k), \bar{y}-y^{k+1} \rangle + \tau_{k+1} \langle A^*(y^{k+1}-y^{k}), x^{k+1}-x^{k+2}\rangle + \tau_{k+1}(\delta_{k+1}+\delta_{k+2}+\varepsilon_{k+1}),$$
which, by the definition of $G(x, y)$, implies (\ref{3.5}) immediately.

\begin{lemma}
Let $\{(z^{k+1}, x^{k+1}, y^{k+1}) : k \geqslant\ 0\}$ be the sequence generated by Algorithm \ref{A1}. For $k \geqslant 0$, it holds that
\begin{equation} \label{3.12}
	\begin{split}
		&\frac{\phi}{\phi-1}\|z^{k+3}-\bar{x} \|_S^2 + \frac{1}{\beta}\| y^{k+1}-\bar{y}\|_T^2 + 2\tau_{k+1}G(x^{k+1}, y^{k+1})\\
		&\leqslant \frac{\phi}{\phi-1}\|z^{k+2}-\bar{x} \|_S^2 + \frac{1}{\beta}\| y^k-\bar{y}\|_T^2 - (1-\frac{\eta}{\lambda_1})\theta_{k+1}\phi\|x^{k+2}-x^{k+1}\|_S^2 \\
		&- \frac{\lambda_2-\eta}{\lambda_2\beta}\|y^{k+1}-y^k\|_T^2  -\theta_{k+1}\phi\|x^{k+1}-z^{k+2}\|_S^2+2\tau_{k+1}(\delta_{k+1}+\delta_{k+2} +\varepsilon_{k+1}).
	\end{split}
\end{equation}
\end{lemma}

\noindent Proof.
First, it is easy to verify from $ \psi=\frac{1+\phi}{\phi^2}$ and $\theta_{k+1}=\frac{\tau_{k+1}}{\tau_{k}} \leqslant \psi$ that
\begin{equation}
1+\frac{1}{\phi} - \phi\theta_{k+1} \geqslant 1+\frac{1}{\phi} - \phi\psi=0.
\end{equation}
Since $\lambda_1$ and $\lambda_2$ are the minimum eigenvalues of $S$ and $T$, respectively, it follows from (\ref{2.10}) that
$$ \|x^{k+2}-x^{k+1} \| \leqslant \frac{1}{\sqrt{\lambda_1}}\|x^{k+2}-x^{k+1} \|_S \quad and \quad \|A^*y^{k+1}-A^*y^k\| \leqslant \frac{1}{\sqrt{\lambda_2}}\|A^*y^{k+1}-A^*y^k\|_T.$$
Using Cauchy-Schwarz inequality, from (\ref{3.4}) we obtain
\begin{equation} \label{3.14}
\begin{split}
	&2\tau_{k+1}\|A^*y^{k+1}-A^*y^k\| \|x^{k+2}-x^{k+1} \| \\ &\leqslant 2\tau_{k+1} \frac{1}{\sqrt{\lambda_1}\sqrt{\lambda_2}}\|x^{k+2}-x^{k+1} \|_S \|A^*y^{k+1}-A^*y^k\|_T \\
	&\leqslant \eta(\frac{\theta_{k+1}\phi}{\lambda_1}\|x^{k+2}-x^{k+1} \|_S^2 + \frac{1}{\lambda_2\beta}\| y^{k+1}-y^k\|_T^2).
\end{split}
\end{equation}
Applying (\ref{2.8}) and Cauchy-Schwarz inequality, from (\ref{3.5}) we have
\begin{equation} \label{3.15}
\begin{split}
	\|x^{k+2}&-\bar{x} \|_S^2 + \frac{1}{\beta}\| y^{k+1}-\bar{y}\|_T^2 + 2\tau_{k+1}G(x^{k+1}, y^{k+1}) \\
	\leqslant &\|z^{k+2}-\bar{x} \|_S^2 + (\phi\theta_{k+1}-1)\|x^{k+2}-z^{k+2} \|_S^2 - \phi\theta_{k+1}\|x^{k+1}-z^{k+2} \|_S^2 \\
	&- \phi\theta_{k+1}\|x^{k+2}-x^{k+1} \|_S^2+\frac{1}{\beta}\| y^k-\bar{y}\|_T^2 - \frac{1}{\beta}\| y^{k+1}-y^k\|_T^2 \\
	&+ 2\tau_{k+1}\|A^*y^{k+1}-A^*y^k\| \|x^{k+1}-x^{k+2}\| +2\tau_{k+1}(\delta_{k+1}+\delta_{k+2} +\varepsilon_{k+1}).
\end{split}
\end{equation}
In view of (\ref{3.1}), we get $ x^{k+2} = \frac{\phi}{\phi-1}z^{k+3} - \frac{1}{\phi-1}z^{k+2}$ and $ z^{k+3} - z^{k+2} = \frac{\phi-1}{\phi}(x^{k+2}-z^{k+2})$.
Thus, from (\ref{2.9}) we deduce
\begin{equation} \label{3.16}
\begin{split}
	\|x^{k+2}-\bar{x} \|_S^2
	&=\|\frac{\phi}{\phi-1}(z^{k+3}-\bar{x})-\frac{1}{\phi-1}(z^{k+2}-\bar{x})\|_S^2\\
	&=\frac{\phi}{\phi-1}\|z^{k+3}-\bar{x} \|_S^2 - \frac{1}{\phi-1}\|z^{k+2}-\bar{x} \|_S^2 + \frac{\phi}{(\phi-1)^2}\|z^{k+3}-z^{k+2} \|_S^2\\
	&=\frac{\phi}{\phi-1}\|z^{k+3}-\bar{x} \|_S^2 - \frac{1}{\phi-1}\|z^{k+2}-\bar{x} \|_S^2 + \frac{1}{\phi}\|x^{k+2}-z^{k+2} \|_S^2.
\end{split}
\end{equation}
Combining (\ref{3.14}) and (\ref{3.16}) with (\ref{3.15}), we obtain
\begin{equation}
\begin{split}
	\frac{\phi}{\phi-1}&\|z^{k+3}-\bar{x} \|_S^2 + \frac{1}{\beta}\| y^{k+1}-\bar{y}\|_T^2 + 2\tau_{k+1}G(x^{k+1}, y^{k+1})\\
	\leqslant &\frac{\phi}{\phi-1}\|z^{k+2}-\bar{x} \|_S^2 -(1-\phi\theta_{k+1}+\frac{1}{\phi})\|x^{k+2}-z^{k+2} \|_S^2+ \frac{1}{\beta}\| y^k-\bar{y}\|_T^2\\
	& - \theta_{k+1}\phi\|x^{k+2}-x^{k+1}\|_S^2 +\eta(\frac{\theta_{k+1}\phi}{\lambda_1}\|x^{k+2}-x^{k+1} \|_S^2 + \frac{1}{\lambda_2\beta}\| y^{k+1}-y^k\|_T^2)\\
	&- \frac{1}{\beta}\|y^{k+1}-y^k\|_T^2  -\theta_{k+1}\phi\|x^{k+1}-z^{k+2}\|_S^2+2\tau_{k+1}(\delta_{k+1}+\delta_{k+2} +\varepsilon_{k+1})\\
	\leqslant &\frac{\phi}{\phi-1}\|z^{k+2}-\bar{x} \|_S^2 + \frac{1}{\beta}\| y^k-\bar{y}\|_T^2 - (1-\frac{\eta}{\lambda_1})\theta_{k+1}\phi\|x^{k+2}-x^{k+1}\|_S^2 \\
	&- \frac{\lambda_2-\eta}{\lambda_2\beta}\|y^{k+1}-y^k\|_T^2 -\theta_{k+1}\phi\|x^{k+1}-z^{k+2}\|_S^2+2\tau_{k+1}(\delta_{k+1}+\delta_{k+2} +\varepsilon_{k+1}).
\end{split}
\end{equation}

In the following, we summarize the convergence result for Algorithm \ref{A1}.
\begin{theorem}
Suppose that $\{(z^{k+1}, x^{k+1}, y^{k+1}) : k \geqslant\ 0\}$ is the sequence generated by Algorithm \ref{A1}. Then, $\{(x^{k+1}, y^{k+1}) : k \geqslant\ 0\}$ is bounded and every limit point of $\{(x^{k+1}, y^{k+1}) : k \geqslant\ 0\}$ is a solution of (\ref{1.1}).
\end{theorem}

\noindent Proof.
Since $\lambda_1,\lambda_2 > \eta$ and $G(x^{k+1},y^{k+1}) \geqslant 0$, (\ref{3.12}) yields
\begin{equation}
\begin{split}
	&\frac{\phi}{\phi-1}\|z^{k+3}-\bar{x} \|_S^2 + \frac{1}{\beta}\| y^{k+1}-\bar{y}\|_T^2 \\
	&\leqslant \frac{\phi}{\phi-1}\|z^{k+2}-\bar{x} \|_S^2 + \frac{1}{\beta}\| y^k-\bar{y}\|_T^2 +2\tau_{k+1}(\delta_{k+1}+\delta_{k+2} +\varepsilon_{k+1}).
\end{split}
\end{equation}
By summing over $k=0,1,2,...,N-1$, we obtain
\begin{equation}
\begin{split}
	&\frac{\phi}{\phi-1}\|z^{N+2}-\bar{x} \|_S^2 + \frac{1}{\beta}\| y^{N}-\bar{y}\|_T^2 \\
	&\leqslant \frac{\phi}{\phi-1}\|z^{2}-\bar{x} \|_S^2 + \frac{1}{\beta}\| y^0-\bar{y}\|_T^2 + 2\sum_{k=0}^{N-1}\tau_{k+1}(\delta_{k+1}+\delta_{k+2} +\varepsilon_{k+1}).
\end{split}
\end{equation}
Since the sequences $\{\delta_{k}\}$ and $\{\varepsilon_{k}\}$ are summable,  and $\{\tau_{k}\}$ is bounded, $\sum_{k=0}^{N-1}\tau_{k+1}(\delta_{k+1}+\delta_{k+2} +\varepsilon_{k+1})$ is bounded. From this we deduce that
$ \frac{\phi}{\phi-1}\|z^{k+1}-\bar{x} \|_S^2 + \frac{1}{\beta}\| y^{k+1}-\bar{y}\|_T^2 $ is bounded. Thus, $\{z^{k+1}\}$ and $\{y^{k+1}\}$ are bounded sequences. Hence, from (\ref{3.1}) we know that $\{x^{k+1}\}$ is also bounded.

Summing up (\ref{3.12}) from $k = 0$ to $N-1$, we get
$$ \frac{\phi}{\phi-1}\|z^{N+2}-\bar{x} \|_S^2 + \frac{1}{\beta}\| y^{N}-\bar{y}\|_T^2 + \phi \sum_{k=0}^{N-1} \theta_{k+1}\|x^{k+1}-z^{k+2}\|_S^2 $$
$$\leqslant \frac{\phi}{\phi-1}\|z^{2}-\bar{x} \|_S^2 + \frac{1}{\beta}\| y^0-\bar{y}\|_T^2 + 2\sum_{k=0}^{N-1}\tau_{k+1}(\delta_{k+1}+\delta_{k+2} +\varepsilon_{k+1}) < \infty.$$
Since $ \theta_{k+1}=\frac{\tau_{k+1}}{\tau_k}>1$,  $\sum_{k=0}^{\infty} \theta_{k+1} = \infty$. Letting $N \to \infty$ in the above inequality and applying the equivalence of $\| \cdot \|_M$ and $\| \cdot \|_2$, where $M$ denotes the symmetric positive definite matrix, we get $\lim\limits_{k \to \infty}\|x^{k+1}-z^{k+2}\| = 0$.

 Similarly, we can deduce that $\lim\limits_{k \to \infty}\|y^{k+1}-y^{k}\| = 0$. Thus, $ (x^{k+1},y^{k+1})$ has at least one limit point $ (x^{\infty},y^{\infty})$ and hence there exists a subsequence $ \{(x^{k_i+1},y^{k_i+1})\}$ such that $ (x^{k_i+1},y^{k_i+1}) \to (x^{\infty},y^{\infty})$ as $i \to \infty$. Similar to (\ref{3.7}) and (\ref{3.11}), for any $(x,y)$, there hold
\begin{equation} \label{3.20}
\tau_{k_i}(f(x^{k_i+1})-f(x)) \leqslant  \langle S(x^{k_i+1}-z^{k_i+1}) + \tau_{k_i}A^*y^{k_i}, x- x^{k_i+1} \rangle + \tau_{k_i}\delta_{k_i+1}.
\end{equation}
and
\begin{equation} \label{3.21}
\tau_{k_i+1}(g(y^{k_i+1})-g(y)) \leqslant \frac{1}{\beta} \langle T(y^{k_i+1} - y^{k_i}) - \beta\tau_{k_i+1}Ax^{k_i+1}, y-y^{k_i+1} \rangle + \tau_{k_i+1}\varepsilon_{k_i+1}.
\end{equation}
 In view of Lemma \ref{L3.1} (ii), we have $\tau_{k_i+1} > \underline{\tau} >0$. Letting $k\rightarrow\infty$  in \eqref{3.20} and \eqref{3.21}, respectively,  and taking into account that $\delta_{k} \to 0$ and $\varepsilon_{k} \to 0$ as $k \to \infty$, we obtain
\begin{equation}
f(x) - f(x^\infty) + \langle A^*y^\infty, x- x^\infty \rangle \geqslant 0 \quad \textrm{and} \quad g(y)-g(y^\infty) - \langle Ax^\infty, y-y^\infty \rangle \geqslant 0.
\end{equation}
This shows that $(x^\infty,y^\infty)$ is a saddle point of (\ref{1.1}). \par
We now establish the convergence rates of Algorithm \ref{A1}.
\begin{theorem}
Let $\{(z^{k+1}, x^{k+1}, y^{k+1}) : k \geqslant\ 0\}$ be the sequence generated by Algorithm \ref{A1}, and $\{{(\bar{x},\bar{y})}\}$ is any saddle point of (\ref{1.1}). Then, for the  ergodic sequence $\{(X^N,Y^N)\}$ given by
\begin{equation} \label{3.23}
	{X}^N = \frac{1}{S^N}\sum_{k=1}^{N}\tau_{k}x^{k}\quad and\quad {Y}^N = \frac{1}{S^N}\sum_{k=1}^{N}\tau_{k}y^{k}\qquad with\quad S^N = \sum_{k=1}^{N}\tau_{k},
\end{equation}
it holds that
\begin{equation} \label{3.24}
	G(X^N,Y^N) \leqslant \frac{c_1}{2\underline{c}N},
\end{equation}
where
$$c_1 = \frac{\phi}{\phi-1}\|z^2 - \bar{x} \|_S^2 + \frac{1}{\beta}\| y^0-\bar{y}\|_T^2 + \sum_{k=0}^{N-1}2\tau_{k+1}(\delta_{k+1}+\delta_{k+2} + \varepsilon_{k+1}). $$
\end{theorem}

\noindent Proof.
Since $\lambda_1, \lambda_2 > \eta $, it follows from (\ref{3.12}) that
\begin{align*}
2\tau_{k+1}G(x^{k+1}, y^{k+1}) \leqslant& \frac{\phi}{\phi-1}(\|z^{k+2}-\bar{x} \|_S^2 - \|z^{k+3}-\bar{x} \|_S^2  ) + \frac{1}{\beta}(\| y^{k}-\bar{y}\|_T^2- \| y^{k+1}-\bar{y}\|_T^2)\\
&+ 2\tau_{k+1}(\delta_{k+1}+\delta_{k+2} + \varepsilon_{k+1}).
\end{align*}
\noindent By taking summation over $k=0,1,2,...,N-1$, we obtain
\begin{equation} \label{3.25}
2\sum_{k=0}^{N-1}\tau_{k+1}G(x^{k+1}, y^{k+1}) \leqslant \frac{\phi}{\phi-1}\|z^2-\bar{x} \|_S^2 + \frac{1}{\beta}\| y^0-\bar{y}\|_T^2 + 2\sum_{k=0}^{N-1}(\tau_{k+1}(\delta_{k+1}+\delta_{k+2} + \varepsilon_{k+1})).
\end{equation}
Since $P(x)$ and $D(y)$ are convex,
	\begin{equation*}
		\begin{cases}
			P(X^N) \leqslant \frac{1}{S^N}\sum_{k=0}^{N-1}\tau_{k+1}P(x^{k+1})  \\
			D(Y^N) \leqslant \frac{1}{S^N}\sum_{k=0}^{N-1}\tau_{k+1}D(y^{k+1}),
		\end{cases}
	\end{equation*}
where $S^N = \sum_{k=0}^{N-1}\tau_{k+1}$. Combining (\ref{3.23}) with (\ref{3.25}), we obtain
\begin{equation}
\begin{split}
	G(&X^N, Y^N)
	= P(X^N)+D(Y^N) \\
	&\leqslant \frac{1}{S^N}\sum_{k=0}^{N-1}\tau_{k+1}(P(x^{k+1})+D(y^{k+1}))=\frac{1}{S^N}\sum_{k=0}^{N-1}\tau_{k+1}G(x^{k+1}, y^{k+1})\\
	& \leqslant \frac{1}{2S^N}(\frac{\phi}{\phi-1}\|z^2-\bar{x} \|_S^2 + \frac{1}{\beta}\| y^0-\bar{y}\|_T^2 + 2\sum_{k=0}^{N-1}(\tau_{k+1}(\delta_{k+1}+\delta_{k+2} + \varepsilon_{k+1}))).
\end{split}
\end{equation}
In view of Lemma \ref{L3.1}(iii), we have $S^N = \sum_{k=1}^{N}\tau_{k} \geqslant \underline{c}N$. Hence, (\ref{3.24}) holds.

\begin{lemma}
{\rm ( \cite{RJ2020})} For $\xi>-1$, let $s^N:=\sum_{k=1}^{N}k^\xi$. Then
$$ s^N = O(N^{1+\xi}). $$
\end{lemma}
Similar to Corallary 3.4 in \cite{RJ2020}, we have the following theorem.
\begin{theorem}
If $\alpha>0$ and $\delta_{k}=O(\frac{1}{k^\alpha})$, $\varepsilon_{k}=O(\frac{1}{k^\alpha})$, we have
	\begin{equation}
      \begin{aligned}
      G(X^N,Y^N)= \left \{
      \begin{array}{ll}
	     O(1/N), \quad \alpha>1,\\
	     O(\ln N/N),  \alpha=1,\\
		 O(N^{-\alpha}), \quad \alpha \in (0,1).
     \end{array} \right.
      \end{aligned}
	\end{equation}
\end{theorem}

\noindent  Proof.
If $\alpha>1$, then the sequences $\{\delta_{k}\}$ and $\{\varepsilon_{k}\}$ are summable. Since the sequence $\{\tau_{k}\}$ is bounded, from (\ref{3.24}) we have
$$ G(X^N,Y^N) =  O(1/N).$$
If $\alpha = 1$, then $\delta_{k} = O(\frac{1}{k}) $ and $\varepsilon_{k} = O(\frac{1}{k})$. In view of the assumption on $\delta_{k}$, for some $r>0$, we have
$$ \delta_{k} \leqslant \frac{r}{k+1}.$$
Thus,
$$ \sum_{k=0}^{N-1} \delta_{k+1} \leqslant \sum_{k=0}^{N-1} \frac{r}{k+1} = c(1+\sum_{k=2}^{N}\frac{1}{k}) \leqslant r(1+\int_{1}^{N}\frac{1}{t}dt) = r(1+\ln N).$$
Hence, from the boundedness of $\{\tau_{k}\}$ we know that $\sum_{k=0}^{N-1} \tau_{k+1}\delta_{k+1} = O(\ln N)$. Similarly, $ \sum_{k=0}^{N-1} \tau_{k+1}\varepsilon_{k+1}$ can also obtain the same result. Therefore, we get
$$ G(X^N,Y^N) =  O(\ln N/N).$$
If $\alpha \in (0,1)$, then $-\alpha > -1$,  from Lemma 3.4 we obtain $\sum_{k=0}^{N-1} \delta_{k+1} = O(N^{1-\alpha})$ and $\sum_{k=0}^{N-1} \varepsilon_{k+1} = O(N^{1-\alpha})$. Thus, we have
$$ G(X^N,Y^N) =  O(N^{-\alpha}). $$

\subsection{Partially strongly convex case}
Now, we consider the acceleration of Algorithm \ref{A1} assuming additionally that $f$ is $\gamma_f$-strongly convex, i.e., it holds for some $\gamma_f>0$ that
\begin{equation} \label{3.28}
f(y) \geqslant f(x) + \langle u,y-x \rangle + \frac{\gamma_f}{2}\|y-x\|^2, \forall u \in \partial f(x), \forall x,y \in \mathbb{R}^{n}.
\end{equation}
When $g$ is strongly convex,the corresponding results can be achieved similarly and is thus omitted. The accelerated version of Algorithm \ref{A1} is summarized in Algorithm \ref{A2}.

\begin{algorithm}[ht]
	1: \quad Choose $x^0 = z^0 \in \mathbb{R}^n $, $y^0 \in \mathbb{R}^m$, $ \beta_0 > 0, \tau_0 >0 , \mu \in (0,1)$ and $\phi \in (\xi,\varphi)$ where $\xi$ is the unique real root of $ \xi^3-\xi-1=0 $. Let $S \in \mathbb{R}^{n \times n}$ and $T \in \mathbb{R}^{m \times m}$ are given symmetric positive definite matrix. $\lambda_1, \lambda_2 > 1 $ are the minimum eigenvalues of $S$ and $T$, respectively, $\Lambda_1$ is the maximum eigenvalues of $S$. Set $ \psi = \frac{1+\phi}{\phi^2}$ and $ k=0.$\\
	2: \quad Compute
	\begin{equation}
		{z}^{k+1} = \frac{\phi-1}{\phi} x^k + \frac{1}{\phi} z^k,
	\end{equation}
	\begin{equation} \label{3.30}
		{x}^{k+1} \approx_2^{\delta_{k+1}} \mathop{\arg\min}\limits_{x \in X} \{ L(x,y^k) + \frac{1}{2\tau_k}\|x-z^{k+1}\|_S^2 \},
	\end{equation}
	\begin{equation} \label{3.31}
		\omega_{k+1} = \frac{\phi - \psi}{\phi\Lambda_1 + \psi\gamma_f\tau_{k}},
	\end{equation}
	\begin{equation} \label{3.32}
		\beta_{k+1} = \beta_k(1+\gamma_f\omega_{k+1}\tau_{k}).
	\end{equation}
	3: \quad Choose any $\tau_{k+1} \in [\tau_{k},\psi\tau_{k}]$ and run\\
	\quad \quad	 3.a: Compute
	\begin{equation} \label{3.33}
		{y}^{k+1} \approx_2^{\varepsilon_{k+1}} \mathop{\arg\max}\limits_{y \in Y} \{ L(x^{k+1},y) - \frac{1}{2\beta_{k+1} \tau_{k+1}}\|y-y^k\|_T^2 \}.
	\end{equation}
	\quad \quad  3.b: Break linesearch if
	\begin{equation}  \label{3.34}
		\sqrt{\beta_{k+1}\tau_{k+1}}\|A^*y^{k+1}-A^*y^k\|_T \leqslant \sqrt{\frac{\phi}{\tau_k}} \|y^{k+1}-y^k\|_T.
\end{equation}
\quad \quad  Otherwise, set $\tau_{k+1}:= \tau_{k+1}\mu$ and go to 3.a.\\
  4: \quad Set $k+1 \leftarrow k+2$ and return to 2.
\caption{Accelerated IP-GRPDAL when $f$ is $\gamma_f$-strongly convex}   \label{A2}
\end{algorithm}

Similar to Lemma \ref{L3.1}, we  have the following results.
\begin{lemma} \label{L3.5}
(i) The linesearch in Algorithm \ref{A2} always terminates.\\
(ii) There exists  constant $c>0$ such that $\beta_k \geqslant ck^2$.
\end{lemma}
\noindent Proof. (i) This conclusion follows  from  Lemma \ref{L3.1}(i) by replacing
$\beta$ with $\beta_k$ and setting $\eta=1$.\\
(ii) The  proof is similar to Lemma 4.1 (ii) in $\cite{C2022}$  and is thus omitted.

Next we will establish the $O(1/{N^2})$ convergence rate of Algorithm \ref{A2}.
\begin{theorem} \label{T3.4}
Let $\{(z^{k+1}, x^{k+1}, y^{k+1}) : k \geqslant\ 0\}$ be the sequence generated by Algorithm \ref{A2}, and $\{{(\bar{x},\bar{y})}\}$ is any saddle point of (\ref{1.1}). Then, we can obtain $\| z^{N+2}-\bar{x} \|_S = O(\frac{1}{N})$ and for the ergodic sequence given by
\begin{equation}
\bar{S}^N = \sum_{k=1}^{N} \beta_k\tau_k,\quad \bar{X}^N=\frac{1}{S^N}\sum_{k=1}^{N}\beta_k\tau_{k}x^k,\quad and \quad \bar{Y}^N=\frac{1}{S^N}\sum_{k=1}^{N}\beta_k\tau_{k}y^k,
\end{equation}
there exists a constant $c_3>0$ such that
\begin{equation}
G(\bar{X}^N,\bar{Y}^N) \leqslant \frac{c_3}{N^2}
\end{equation}
\end{theorem}

\noindent Proof.
Since $f$ is strongly convex, it follows from (\ref{3.30}) and Definition \ref{D2.3} that
\begin{equation}
f(x)-f(x^{k+1})+ \frac{1}{\tau_{k}} \langle {S(x^{k+1}-z^{k+1}) + A^*y^k, x-x^{k+1}} \rangle - \frac{\gamma_f}{2}\| x-x^{k+1} \|^2+ \delta_{k+1} \geqslant 0.
\end{equation}
Since $\Lambda_1$ is the maximum eigenvalue of $S$, from (\ref{2.10}) we get $\| x-x^{k+1} \|^2 \geqslant \frac{1}{\Lambda_1}\| x-x^{k+1} \|_S^2$. Similar to (\ref{3.7})-(\ref{3.10}), we obtain
\begin{equation}  \label{3.38}
	\begin{split}
\tau_{k+1}&(f(x^{k+2})-f(\bar{x})+\frac{\gamma_f}{2\Lambda_1}\|\bar{x}-x^{k+2} \|_S^2)) \\
&\leqslant \langle  S(x^{k+2} - z^{k+2}) + \tau_{k+1}A^*y^{k+1}, \bar{x}-x^{k+2} \rangle + \tau_{k+1}\delta_{k+2},
\end{split}
\end{equation}
\begin{equation} \label{3.39}
		\begin{split}
 \tau_{k+1}&(f(x^{k+1})-f(x^{k+2})+\frac{\gamma_f}{2\Lambda_1}\|x^{k+2}-x^{k+1}\|_S^2)) \\
 &\leqslant \langle \phi\theta_{k+1} S(x^{k+1} - z^{k+2}) + \tau_{k+1}A^*y^{k}, x^{k+2}-x^{k+1} \rangle + \tau_{k+1}\delta_{k+1}.
\end{split}
\end{equation}
From (\ref{3.33}), for $\bar{y} \in Y$ we have
\begin{equation}  \label{3.40}
\tau_{k+1}(g(y^{k+1})-g(\bar{y})) \leqslant \frac{1}{\beta_{k+1}}\langle T(y^{k+1} - y^k) x - \beta_{k+1}\tau_{k+1}Ax^{k+1}, \bar{y}-y^{k+1} \rangle + \tau_{k+1}\varepsilon_{k+1}.
\end{equation}
Then, by adding (\ref{3.38}), (\ref{3.39}) and (\ref{3.40}) and using similar arguments as in Lemma \ref{L3.2}, we deduce
\begin{equation}  \label{3.41}
\begin{split}
\tau_{k+1}G&(x^{k+1},y^{k+1}) \leqslant  \langle S(x^{k+2} - z^{k+2}), \bar{x}-x^{k+2} \rangle+ \phi\theta_{k+1} \langle S(x^{k+1} - z^{k+2}), x^{k+2}-x^{k+1} \rangle \\
&+ \frac{1}{\beta_{k+1}}\langle T(y^{k+1} - y^k), \bar{y}-y^{k+1} \rangle + \tau_{k+1}\langle A^*(y^{k+1}-y^{k}), x^{k+1}-x^{k+2}\rangle\\
&-\frac{\gamma_f\tau_{k+1}}{2\Lambda_1}\|x^{k+2}-x^{k+1}\|_S^2-\frac{\gamma_f\tau_{k+1}}{2\Lambda_1}\|\bar{x}-x^{k+2}\|_S^2+\tau_{k+1}(\delta_{k+1}+\delta_{k+2}+\varepsilon_{k+1}).
\end{split}
\end{equation}
Using (\ref{2.8}) and Cauchy-Schwarz inequality, from (\ref{3.41}) we get
\begin{equation}  \label{3.42}
\begin{split}
(1+&\frac{\gamma_f \tau_{k+1}}{\Lambda_1})\|\bar{x}-x^{k+2}\|_S^2 + \frac{1}{\beta_{k+1}}\|\bar{y}-y^{k+1}\|_T^2 + 2\tau_{k+1}G(x^{k+1},y^{k+1}) \\
 \leqslant &\|z^{k+2}-\bar{x}\|_S^2 + \frac{1}{\beta_{k+1}}\|y^k-\bar{y}\|_T^2 - \frac{1}{\beta_{k+1}}\|y^k-y^{k+1}\|_T^2 \\
& + (\phi\theta_{k+1}-1)\|z^{k+2}-x^{k+2}\|_S^2 - \phi\theta_{k+1}\|z^{k+2}-x^{k+1}\|_S^2- \phi\theta_{k+1}\|x^{k+2}-x^{k+1}\|_S^2 \\
&+ 2\tau_{k+1}\|A^*(y^{k+1}-y^k)\|\| x^{k+1}-x^{k+2}\| + 2\tau_{k+1}(\delta_{k+1}+\delta_{k+2}+\varepsilon_{k+1}).
\end{split}
\end{equation}
Combining (\ref{3.16}) with (\ref{3.42}), we have
\begin{equation} \label{3.43}
\begin{split}
(1+&\frac{\gamma_f \tau_{k+1}}{\Lambda_1})\frac{\phi}{\phi-1}\|\bar{x}-z^{k+3}\|_S^2 + \frac{1}{\beta_{k+1}}\|\bar{y}-y^{k+1}\|_T^2 + 2\tau_{k+1}G(x^{k+1},y^{k+1}) \\
\leqslant & \frac{(\gamma_f\tau_{k+1} / \Lambda_1)+\phi}{\phi-1}\|z^{k+2}-\bar{x}\|_S^2 + \frac{1}{\beta_{k+1}}\|y^k-\bar{y}\|_T^2 - \frac{1}{\beta_{k+1}}\|y^k-y^{k+1}\|_T^2 \\
& + (\phi\theta_{k+1}-1 -\frac{1+(\gamma_f\tau_{k+1} / \Lambda_1)}{\phi})\|z^{k+2}-x^{k+2}\|_S^2 - \phi\theta_{k+1}\|z^{k+2}-x^{k+1}\|_S^2  \\
&- \phi\theta_{k+1}\|x^{k+2}-x^{k+1}\|_S^2+ 2\tau_{k+1}\|A^*(y^{k+1}-y^k)\|\| x^{k+1}-x^{k+2}\| \\
&+ 2\tau_{k+1}(\delta_{k+1}+\delta_{k+2}+\varepsilon_{k+1}).
\end{split}
\end{equation}
Thus, it follows from (\ref{2.10}), (\ref{3.34}) and Cauchy-Schwarz inequality that
\begin{equation} \label{3.44}
2\tau_{k+1}\|A^*y^{k+1}-A^*y^k\| \|x^{k+1}-x^{k+2} \| \leqslant \frac{\theta_{k+1}\phi}{\lambda_1}\|x^{k+2}-x^{k+1} \|_S^2 + \frac{1}{\lambda_2\beta_{k+1}}\| y^{k+1}-y^k\|_T^2.
\end{equation}
Since $\psi=\frac{1+\phi}{\phi^2}$ and $ \theta_{k+1} \leqslant \psi$, we have $\phi\theta_{k+1}-1 -\frac{1+(\gamma_f\tau_{k+1} / \Lambda_1)}{\phi} \leqslant \phi\psi-1-\frac{1}{\phi}-\frac{(\gamma_f\tau_{k+1} / \Lambda_1)}{\phi}=-\frac{\gamma_f\tau_{k+1}}{\Lambda_1\phi}$. Therefore, substituting (\ref{3.44}) into (\ref{3.43}), we obtain
\begin{equation} \label{3.45}
\begin{split}
(1+&\frac{\gamma_f \tau_{k+1}}{\Lambda_1})\frac{\phi}{\phi-1}\|\bar{x}-z^{k+3}\|_S^2 + \frac{1}{\beta_{k+1}}\|\bar{y}-y^{k+1}\|_T^2 + 2\tau_{k+1}G(x^{k+1},y^{k+1}) \\
\leqslant & \frac{(\gamma_f\tau_{k+1}/\Lambda_1)+\phi}{\phi-1}\|z^{k+2}-\bar{x}\|_S^2  + \frac{1}{\beta_{k+1}}\|y^k-\bar{y}\|_T^2 - \frac{\gamma_f\tau_{k+1}}{\phi\Lambda_1}\|z^{k+2}-x^{k+2}\|_S^2 \\
&-(1-\frac{1}{\lambda_1})\theta_{k+1}\phi\|x^{k+2}-x^{k+1}\|_S^2-\frac{\lambda_2-1}{\lambda_2\beta_{k+1}}\|y^{k+1}-y^k\|_T^2+ 2\tau_{k+1}(\delta_{k+1}+\delta_{k+2}+\varepsilon_{k+1})\\
\leqslant& \frac{(\gamma_f\tau_{k+1}/\Lambda_1)+\phi}{\phi-1}\|z^{k+2}-\bar{x}\|_S^2  + \frac{1}{\beta_{k+1}}\|y^k-\bar{y}\|_T^2 - \frac{\gamma_f\tau_{k+1}}{\phi\Lambda_1}\|z^{k+2}-x^{k+2}\|_S^2 \\
&+ 2\tau_{k+1}(\delta_{k+1}+\delta_{k+2}+\varepsilon_{k+1}).
\end{split}
\end{equation}
Note that $ (1+\frac{\gamma_f \tau_{k+1}}{\Lambda_1})\frac{\phi}{\phi-1} = \frac{\phi(1+(\gamma_f\tau_{k+1}/\Lambda_1))}{\phi+(\gamma_f\tau_{k+2}/\Lambda_1)}\frac{\phi+(\gamma_f\tau_{k+2}/\Lambda_1)}{\phi-1}.$ Thus, it follows from $ \tau_{k+2} \leqslant \psi\tau_{k+1}$ and $ \omega_{k+2} = \frac{\phi - \psi}{\phi\Lambda_1 + \psi\gamma_f\tau_{k+1}} $ that
\begin{equation}
\frac{\phi(1+(\gamma_f\tau_{k+1}/\Lambda_1))}{\phi+(\gamma_f\tau_{k+2}/\Lambda_1)} \geqslant \frac{\phi(1+(\gamma_f\tau_{k+1}/\Lambda_1))}{\phi+(\gamma_f\psi\tau_{k+1}/\Lambda_1)} = 1+\frac{(\phi-\psi)\gamma_f\tau_{k+1}}{\phi\Lambda_1+\gamma_f\psi\tau_{k+1}} = 1+\omega_{k+2}\gamma_f\tau_{k+1}.
\end{equation}
Since $ \beta_{k+2} = \beta_{k+1}(1+\gamma_f\omega_{k+2}\tau_{k+1})$,
\begin{equation}  \label{3.47}
\begin{split}
(1+\frac{\gamma_f \tau_{k+1}}{\Lambda_1})\frac{\phi}{\phi-1}
&=\frac{\phi(1+(\gamma_f\tau_{k+1}/\Lambda_1))}{\phi+(\gamma_f\tau_{k+2}/\Lambda_1)}\frac{\phi+(\gamma_f\tau_{k+2}/\Lambda_1)}{\phi-1}\\
&\geqslant (1+\omega_{k+2}\gamma_f\tau_{k+1})\frac{\phi+(\gamma_f\tau_{k+2}/\Lambda_1)}{\phi-1} \\
&= \frac{\beta_{k+2}}{\beta_{k+1}}\frac{\phi+(\gamma_f\tau_{k+2}/\Lambda_1)}{\phi-1}.
\end{split}
\end{equation}
Substituting (\ref{3.47}) into (\ref{3.45}) and multiplying the resulting inequality by $\frac{1}{2}\beta_{k+1}$, we have
\begin{equation}  \label{3.48}
\begin{split}
\frac{\beta_{k+2}}{2}&\frac{\phi+(\gamma_f\tau_{k+2}/\Lambda_1)}{\phi-1} \|\bar{x}-z^{k+3}\|_S^2 + \frac{1}{2}\|\bar{y}-y^{k+1}\|_T^2 + \beta_{k+1}\tau_{k+1}G(x^{k+1},y^{k+1})\\
\leqslant & \frac{\beta_{k+1}}{2}\frac{(\gamma_f\tau_{k+1}/\Lambda_1)+\phi}{\phi-1}\|z^{k+2}-\bar{x}\|_S^2 + \frac{1}{2}\|y^k-\bar{y}\|_T^2  \\
&- \frac{\beta_{k+1}}{2}\frac{\gamma_f\tau_{k+1}}{\Lambda_1\phi}\|z^{k+2}-x^{k+2}\|_S^2+\beta_{k+1}\tau_{k+1}(\delta_{k+1}+\delta_{k+2}+\varepsilon_{k+1}).
\end{split}
\end{equation}
Define $A_{k+1} := \frac{(\gamma_f\tau_{k+1}/\Lambda_1)+\phi}{2(\phi-1)}\|z^{k+2}-\bar{x}\|_S^2 + \frac{1}{2\beta_{k+1}}\|y^k-\bar{y}\|_T^2 $, and hence  (\ref{3.48}) yields
\begin{align*}
	\beta_{k+2}A_{k+2} + \beta_{k+1}\tau_{k+1}G(x^{k+1},y^{k+1}) \leqslant &\beta_{k+1}A_{k+1} - \frac{\beta_{k+1}\gamma_f\tau_{k+1}}{2\phi\Lambda_1}\|z^{k+2}-x^{k+2}\|_S^2\\
    &+\beta_{k+1}\tau_{k+1}(\delta_{k+1}+\delta_{k+2}+\varepsilon_{k+1}).
\end{align*}
By summing over $k=0,1,2,...,N-1$ in the above inequality, we obtain
\begin{equation}  \label{3.49}
\begin{split}
\beta_{N+1}A_{N+1} & + \sum_{k=0}^{N-1} \beta_{k+1}\tau_{k+1}G(x^{k+1}, y^{k+1}) + \sum_{k=0}^{N-1}\frac{\beta_{k+1}\gamma_f\tau_{k+1}}{2\phi\Lambda_1}\|z^{k+2}-x^{k+2}\|_S^2 \\
&\leqslant \beta_{1} A_1 + \sum_{k=0}^{N-1}\beta_{k+1}\tau_{k+1}(\delta_{k+1}+\delta_{k+2}+\varepsilon_{k+1}).
\end{split}
\end{equation}
In view of the convexity of $G(x,y)$, from (\ref{3.49}), we have
\begin{equation}
\begin{split}
 G(\bar{X}^N,\bar{Y}^N) &\leqslant \frac{1}{\bar{S}^N} \sum_{k=0}^{N-1}\beta_{k+1}\tau_{k+1}G(x^{k+1},y^{k+1}) \\
& \leqslant \frac{1}{\bar{S}^N}(\beta_{1}A_1 + \sum_{k=0}^{N-1}\beta_{k+1}\tau_{k+1}(\delta_{k+1}+\delta_{k+2}+\varepsilon_{k+1})).
\end{split}
\end{equation}
According to the definition of $A_{k+1}$ and (\ref{3.49}), we get
$$A_{N+1} = \frac{(\gamma_f\tau_{N+1}/\Lambda_1)+\phi}{2(\phi-1)}\|z^{N+2}-\bar{x}\|_S^2 + \frac{1}{2\beta_{N+1}}\|y^N-\bar{y}\|_T^2$$ $$ \leqslant \frac{\beta_{1} A_1 + \sum_{k=0}^{N-1}\beta_{k+1}\tau_{k+1}(\delta_{k+1}+\delta_{k+2}+\varepsilon_{k+1})}{\beta_{N+1}}.$$
Thus,
\begin{equation} \label{3.51}
\| z^{N+2}-\bar{x} \|_S^2 \leqslant \frac{2(\phi-1)}{(\gamma_f\tau_{N+1}/\Lambda_1)+\phi} A_{N+1} \leqslant \frac{2(\beta_{1} A_1 + \sum_{k=0}^{N-1}\beta_{k+1}\tau_{k+1}(\delta_{k+1}+\delta_{k+2}+\varepsilon_{k+1})}{\beta_{N+1}}.
\end{equation}
From Lemma \ref{L3.5}(ii), we know that there exists $c>0$ such that $ \beta_k \geqslant ck^2 $ for all $k \geqslant 1$, and hence (\ref{3.51}) implies $\| z^{N+2}-\bar{x} \|_S \leqslant \frac{c_2}{N}$ with $$c_2=\sqrt{\frac{2(\beta_{1}A_1 + \sum_{k=0}^{N-1}\beta_{k+1}\tau_{k+1}(\delta_{k+1}+\delta_{k+2}+\varepsilon_{k+1}))}{c}}>0.$$
Since $\beta_{k+1} = \beta_k(1+\gamma_f\omega_{k+1}\tau_{k})$ and $\omega_{k+1}=\frac{\phi - \psi}{\phi\Lambda_1 + \psi\gamma_f\tau_{k}}<1$,
\begin{equation}
\beta_k\tau_k = \frac{\beta_{k+1}-\beta_k}{\gamma_f\omega_{k+1}} \geqslant \frac{\beta_{k+1}-\beta_k}{\gamma_f}.
\end{equation}
Since $ \bar{S}^N = \sum_{k=1}^{N} \beta_k\tau_k \geqslant \frac{\beta_{N+1}-\beta_{1}}{\gamma_f}$ and $\beta_k \geqslant ck^2,$ $\bar{S}^N =O(N^2)$. This means that $G(\bar{X}^N,\bar{Y}^N) \leqslant \frac{c_3}{N^2}$ for some constant $c_3>0$. This completes the proof.

\subsection{Completely strongly convex case}
We further assume that $f$ is $\gamma_f$-strongly convex and $g$ is $\gamma_g$-strongly convex. In this setting, Algorithm \ref{A3} can be accelerated to linear convergence by properly selecting parameters $\tau_{k}$, $\delta_{k}$ and $\varepsilon_{k}$.\\
\begin{algorithm}[H]  \label{A3}
\caption{Accelerated IP-GRPDAL when $f$ and $g$ are strongly convex}
1: \quad Choose $x^0 = z^0 \in \mathbb{R}^n $, $y^0 \in \mathbb{R}^m$,  $\beta > 0, \tau_0 >0, \mu \in (0,1)$ and $\phi \in (\xi,\varphi)$  where $\xi$ is the unique real root of $ \xi^3-\xi-1=0 $. Let $S \in \mathbb{R}^{n \times n}$ and $T \in \mathbb{R}^{m \times m}$ are given symmetric positive definite matrix. Set $ \psi = \frac{1+\phi}{\phi^2}$ and $ k=0.$\\
2: \quad Compute
$$ {z}^{k+1} = \frac{\phi-1}{\phi} x^k + \frac{1}{\phi} z^k,$$
$$ 	{x}^{k+1} \approx_2^{\delta_{k+1}} \mathop{\arg\min}\limits_{x \in X} \{ L(x,y^k) + \frac{1}{2\tau_k}\|x-z^{k+1}\|_S^2 \}.$$\\
3: \quad Choose any $\tau_{k+1} \in [\tau_{k},\psi\tau_{k}]$ and run\\
\quad \quad	 3.a: Compute
\begin{equation}
	{y}^{k+1} \approx_2^{\varepsilon_{k+1}} \mathop{\arg\max}\limits_{y \in Y} \{ L(x^{k+1},y) - \frac{1}{2\beta \tau_{k+1}}\|y-y^k\|_T^2 \}.
\end{equation}
\quad \quad  3.b: Break linesearch if
\begin{equation}
	\sqrt{\beta\tau_{k+1}}\|A^*y^{k+1}-A^*y^k\|_T \leqslant \sqrt{\frac{\phi}{\tau_k}} \|y^{k+1}-y^k\|_T.
\end{equation}
\quad \quad  Otherwise, set $\tau_{k+1}:= \tau_{k+1}\mu$ and go to 3.a.\\
4: \quad  Set $k+1 \leftarrow k+2$ and return to 2.
\end{algorithm}

Now, we summarize the linear convergence rate of Algorithm \ref{A3} in the following theorem.
\begin{theorem}
Suppose that $\Lambda_1,\Lambda_2 >1$ are the maximum eigenvalues of $S$ and $T$, respectively, $ \delta_{k}=\varepsilon_{k}=O(q^k)$ with $q \in (0,1)$, and $ \tau_{k}= \tau$ such that $ 1+\frac{\gamma_f\tau}{\Lambda_1}=1+\frac{\beta\gamma_g\tau}{\Lambda_2}=\frac{1}{\rho}$. Let $\{(z^{k+1}, x^{k+1}, y^{k+1}) : k \geqslant\ 0\}$ be the sequence generated by Algorithm \ref{A3} and $\{(\bar{x},\bar{y})\}$ is any saddle point of (\ref{1.1}).  Then, for the ergodic sequence $\{(\widetilde{X}^N,\widetilde{Y}^N)\}$ given by
\begin{equation} \label{3.55}
\widetilde{X}^N= \frac{1}{\widetilde{S}^N}\sum_{k=0}^{N-1}\frac{1}{\rho ^k}x^{k+1} \quad and \quad \widetilde{Y}^N= \frac{1}{\widetilde{S}^N}\sum_{k=0}^{N-1}\frac{1}{\rho ^k}y^{k+1} \quad with \quad \widetilde{S}^N= \sum_{k=0}^{N-1} \frac{1}{\rho^k},
\end{equation}
it holds that
	\begin{equation*}
      \begin{aligned}
      G(\widetilde{X}^N,\widetilde{Y}^N)+\frac{1}{2\tau}\|\bar{x}-z^{N+2}\|_S^2= \left \{
      \begin{array}{ll}
	     O(\rho^N), \quad \rho>q,\\
	    O(N\rho^N),  \rho=q,\\
		O(q^N), \quad \rho<q.
     \end{array} \right.
      \end{aligned}
	\end{equation*}
\end{theorem}
\noindent  Proof.
Similar to the proof of (\ref{3.41}) in Theorem \ref{T3.4}, one can deduce that
\begin{equation}
\begin{split}
\tau G( & x^{k+1},y^{k+1}) \leqslant \langle S( x^{k+2} - z^{k+2}), \bar{x}-x^{k+2} \rangle+ \phi\theta_{k+1}\langle S (x^{k+1} - z^{k+2}), x^{k+2}-x^{k+1} \rangle \\
&+ \frac{1}{\beta}\langle T(y^{k+1} - y^k), \bar{y}-y^{k+1} \rangle + \tau\langle A^*(y^{k+1}-y^{k}), x^{k+1}-x^{k+2}\rangle -\frac{\gamma_f\tau}{2\Lambda_1}\|x^{k+2}-\bar{x}\|_S^2\\
&-\frac{\gamma_f\tau}{2\Lambda_1}\|x^{k+2}-x^{k+1}\|_S^2-\frac{\gamma_g\tau}{2\Lambda_2}\|y^{k+1}-\bar{y}\|_T^2+\tau(\delta_{k+1}+\delta_{k+2}+\varepsilon_{k+1}).
\end{split}
\end{equation}
Applying the same arguments as in (\ref{3.42})-(\ref{3.45}), we have
\begin{equation} \label{3.57}
\begin{split}
(1&+\frac{\gamma_f\tau}{\Lambda_1})\frac{\phi}{\phi-1}\|z^{k+3}-\bar{x}\|_S^2 + (\frac{1}{\beta}+\frac{\gamma_g\tau}{\Lambda_2})\|y^{k+1}-\bar{y}\|_T^2 + 2\tau G(x^{k+1},y^{k+1}) \\
\leqslant&  \frac{(\gamma_f\tau/\Lambda_1)+\phi}{\phi-1}\|z^{k+2}-\bar{x}\|_S^2  + \frac{1}{\beta}\|y^k-\bar{y}\|_T^2 -\frac{\gamma_f\tau}{\phi\Lambda_1}\|z^{k+2}-x^{k+2}\|_S^2\\
&+ 2\tau(\delta_{k+1}+\delta_{k+2}+\varepsilon_{k+1}).
\end{split}
\end{equation}
 Since $1+\frac{\gamma_f\tau}{\Lambda_1}=1+\frac{\beta\gamma_g\tau}{\Lambda_2}=\frac{1}{\rho}$, from (\ref{3.47}) we get
$$ (1+\frac{\gamma_f\tau}{\Lambda_1})\frac{\phi}{\phi-1} \geqslant \frac{\phi(1+(\gamma_f\tau/\Lambda_1))}{\phi+(\gamma_f\tau/\Lambda_1)}\frac{\phi+(\gamma_f\tau/\Lambda_1)}{\phi-1} \geqslant \frac{1}{\rho}\frac{\phi+(\gamma_f\tau/\Lambda_1)}{\phi-1} \quad and \quad \frac{1}{\beta}+\frac{\gamma_g\tau}{\Lambda_2} = \frac{1}{\rho\beta}.$$
Thus, the inequality (\ref{3.57}) implies that
\begin{equation} \label{3.58}
\begin{split}
\frac{1}{\rho}&\frac{\phi+(\gamma_f\tau/\Lambda_1)}{\phi-1}\|z^{k+3}-\bar{x}\|_S^2 + \frac{1}{\rho\beta}\|y^{k+1}-\bar{y}\|_T^2 + 2\tau G(x^{k+1},y^{k+1}) \\
&\leqslant \frac{\phi+(\gamma_f\tau/\Lambda_1)}{\phi-1}\|z^{k+2}-\bar{x}\|_S^2  + \frac{1}{\beta}\|y^k-\bar{y}\|_T^2 + 2\tau (\delta_{k+1}+\delta_{k+2}+\varepsilon_{k+1}).
\end{split}
\end{equation}
Multiplying (\ref{3.58}) by $\rho^{-k}$ and summing the resulting inequality from $k = 0$ to $N-1$, we get
\begin{equation} \label{3.59}
\begin{split}
\frac{1}{\rho^N}&\frac{\phi+(\gamma_f\tau/\Lambda_1)}{\phi-1}\|z^{N+2}-\bar{x}\|_S^2 + \frac{1}{\rho^N\beta}\|y^{N}-\bar{y}\|_T^2 + 2\tau\sum_{k=0}^{N-1}\frac{1}{\rho^k}G(x^{k+1},y^{k+1}) \\
&\leqslant \frac{\phi+(\gamma_f\tau/\Lambda_1)}{\phi-1}\|z^{2}-\bar{x}\|_S^2  + \frac{1}{\beta}\|y^0-\bar{y}\|_T^2 + 2\tau\sum_{k=0}^{N-1}\frac{1}{\rho^k} (\delta_{k+1}+\delta_{k+2}+\varepsilon_{k+1}).
\end{split}
\end{equation}
Since $P(x)$ and $P(y)$ are convex, from (\ref{3.55}) we deduce
$$ \sum_{k=0}^{N-1}\frac{1}{\rho^k}G(x^{k+1},y^{k+1}) =\sum_{k=0}^{N-1}\frac{1}{\rho^k}(P(x^{k+1})+D(y^{k+1})) \geqslant  \widetilde{S}^N(P(\widetilde{X}^N)+D(\widetilde{Y}^N))=\widetilde{S}^NG(\widetilde{X}^N,\widetilde{Y}^N).$$
Therefore, from (\ref{3.59}) we have
\begin{equation}
\begin{split}
\widetilde{S}^N& G(\widetilde{X}^N,\widetilde{Y}^N) + \frac{\phi+(\gamma_f\tau/\Lambda_1)}{2\tau\rho^N(\phi-1)}\|\bar{x}-z^{N+2}\|_S^2\\
& \leqslant \frac{(\gamma_f\tau_{1}/\Lambda_1)+\phi}{2\tau(\phi-1)}\|z^{2}-\bar{x}\|_S^2 + \frac{1}{2\beta\tau}\|y^0-\bar{y}\|_T^2 + \sum_{k=0}^{N-1}\frac{1}{\rho^k}(\delta_{k+1}+\delta_{k+2}+\varepsilon_{k+1}).
\end{split}
\end{equation}
Note that $\delta_{k}=\varepsilon_{k}=O(q^k)$, $\frac{1}{\widetilde{S}^N}=1/( \sum_{k=0}^{N-1} \frac{1}{\rho^k}) = O(\rho^N)$ with $q,\rho \in (0,1)$, and
$$ \rho^N \sum_{k=0}^{N-1}\frac{q^k}{\rho^k} = \rho^N \sum_{k=0}^{N-1}(\frac{q}{\rho})^k = (q^N-\rho^N)\frac{\rho}{q-\rho}.$$
Thus,
	\begin{equation*}
      \begin{aligned}
      \frac{1}{\widetilde{S}^N}\sum_{k=0}^{N-1}\frac{1}{\rho^k}(\delta_{k+1}+\delta_{k+2}+\varepsilon_{k+1}) = \left \{
      \begin{array}{ll}
	     O(\rho^N), \quad \rho>q,\\
	    O(N\rho^N),  \rho=q,\\
		O(q^N), \quad \rho<q.
     \end{array} \right.
      \end{aligned}
	\end{equation*}
Since $G(\widetilde{X}^N,\widetilde{Y}^N) \geqslant 0$ and $\frac{\phi+(\gamma_f\tau/\Lambda_1)}{\rho^N(\phi-1)} >\frac{\phi+(\gamma_f\tau/\Lambda_1)}{\phi-1} >1 $,
	\begin{equation*}
      \begin{aligned}
     G(\widetilde{X}^N,\widetilde{Y}^N) = \left \{
      \begin{array}{ll}
	     O(\rho^N), \quad \rho>q,\\
	    O(N\rho^N),  \rho=q,\\
		O(q^N), \quad \rho<q,
     \end{array} \right.
      \end{aligned}
           \begin{aligned}
    \frac{1}{2\tau}\|\bar{x}-z^{N+2}\|_S^2 = \left \{
      \begin{array}{ll}
	     O(\rho^N), \quad \rho>q,\\
	    O(N\rho^N),  \rho=q,\\
		O(q^N), \quad \rho<q.
     \end{array} \right.
      \end{aligned}
	\end{equation*}
This completes the proof.

\section{Numerical experiments}\label{sec:exp}
In this section, we apply the proposed Algorithm \ref{A1} to solve sparse recovery and image deblurring problems. We compare IP-GRPDAL method with  Algorithm 1 in \cite{CP2011}(denoted by PDA), Algorithm 1 in \cite{FH2022}(denoted by IPDA), PDAL method in \cite{YM2018}, and GRPDAL method in \cite{C2022}. All codes were written in MATLAB R2015a on a PC (with CPU Intel i5-5200U). For simplicity we set $ S =
\left(
 \begin{array}{cc}
 \frac{1}{r_1}I_1 & 0 \\
0 & \frac{1}{r_2}I_2
\end{array}
\right) $, $ T =
\left(
 \begin{array}{cc}
 \frac{1}{s_1}I_1 & 0 \\
0 & \frac{1}{s_2}I_2
\end{array}
\right) $ , $\delta_{k+1} = 0$ and $\varepsilon_{k+1} = O(1/(k+1)^\alpha)$ in Algorithm \ref{A1}.

\subsection{$l_1$ regularized analysis sparse recovery problem}
We study the following problem:
\begin{equation} \label{4.1}
\varPhi(x) = \mathop{ \textrm{min}}\limits_{x }  \frac{1}{2} \|Ax-b\|^2 + \zeta \|x\|_1
\end{equation}
where $A \in \mathbb{R}^{m \times n}$, $x \in \mathbb{R}^n$, $b \in \mathbb{R}^m$, and $\zeta>0$ is a regularization parameter. Let $f(x) = \zeta \|x\|_1$, $g^*(y) = \frac{1}{2}\|y\|^2 + \langle b,y \rangle$, We can rewrite (\ref{4.1}) as
\begin{equation}
\mathop{ \textrm{min}}\limits_{x } \mathop{ \textrm{max}}\limits_{y} \langle Ax,y \rangle + f(x) - g^*(y)
\end{equation}
We use Proximal Gradient Method to solve subproblems in these algorithms, and the parameters are set as follows:\\
(i) $A := \frac{1}{\sqrt{n}}\textrm{randn}(n,p)$;\\
(ii) $\omega \in R^n$ is a random vector, where $s$ random coordinates are drawn from the uniform distribution in $[-10,10]$ and the rest are zeros;\\
(iii) The observed value $b$ is generated by $b = A\omega + N(0,0.1)$, the regularization parameter $ \zeta $ was set to be $0.1$;\\
(iv) We set $\tau=\frac{1}{10\|A\|}, \sigma=\frac{10}{\|A\|}$ for the PDA method. For  PDAL, GRPDAL and IP-GRPDAL methods, we set $\tau_0=\|y_{-1}-y_0\|/(\sqrt{\beta}\|A^*y_{-1}-A^*y_0\|),\beta=100, \phi=1.618, \mu=0.7$ and $ \eta=0.99.$  At the same time, we set  $s_1 = 2, r_1 = \frac{0.99}{s_1}, s_2 = 1, r_2 = \frac{0.99}{s_2}$ and $\alpha =2$ for the IP-GRPDAL method as those for the IPDA method. The initial points are $x^0=(0,...,0)$ and $y^0=Ax^0+b.$ \par
In this experiment, we terminate all the algorithms when $\varPhi(x^k)-\varPhi^*<10^{-10}$, where an approximation of the optimal value $\varPhi^*$ is obtained by running the algorithms for 5000 iterations.  \par
In order to investigate the stability and efficiency of our method, we test three groups of problems with different pairs $(n, p,s)$ and run the tests 10 times for each group. The average numerical performances including the CPU time (Time, in seconds), the number of iterations (Iter) and the number of extra linesearch trial steps (LS) of PDAL, GRPDAL and IP-GRPDAL methods are reported in Table \ref{tab:table-np}.
\begin{table}[ht]
\centering
\caption{Numerical results of tested algorithms with random tight frames} \label{tab:table-np}
\begin{tabular}{c@{\hspace{7pt}}c@{\hspace{7pt}}cccccccc}
\toprule[1.5pt]
\multirow{3}*[0.9ex]{$n$} & \multirow{3}*[0.9ex]{$p$}& \multirow{3}*[0.9ex]{$s$} & \multicolumn{2}{c}{PDA} & \multicolumn{2}{c}{IPDA} & \multicolumn{3}{c}{PDAL}\\
\cmidrule(r){4-5} \cmidrule(r){6-7} \cmidrule(r){8-10}
&\textbf{ } &\textbf{ } & Time & Iter & Time & Iter & Time & Iter & LS\\
\midrule[1pt]
100 &100 &10 &62.23 &58842 &13.68 &51229 &29.27 &28316 &27918\\
500 &800 &50 &138.89 &$\diagup$ &60.64 &92863 &79.56 &51500 & 50276\\
1000 &2000 &100 &427.67 &$\diagup$  &391.08 &184676 &206.37 &108251 &106938\\
\midrule[1pt]
\multirow{3}*[0.9ex]{$n$} & \multirow{3}*[0.9ex]{$p$}& \multirow{3}*[0.9ex]{$s$} & \multicolumn{3}{c}{GRPDAL} & \multicolumn{4}{c}{\qquad IP-GRPDAL}\\
\cmidrule(r){4-6} \cmidrule(r){8-10}
&\textbf{ } &\textbf{ } & Time & Iter & LS &\textbf{ } & Time & Iter &LS \\
100 & 100 &10 &18.75 &19390 &5976 &\textbf{ } &9.49 &18292 &4647\\
500 &800 &50 &58.23 &42407 &13980 &\textbf{ } &53.01 &38594 &12260\\
1000 &2000 &100 &175.49 &$\diagup$  &72619 &\textbf{ } &113.30 &$\diagup$  &58627\\
\bottomrule[1.5pt]
\end{tabular}
\end{table}

\begin{figure}[!ht]
\centering
\begin{minipage}[t]{0.33\textwidth}
\centering
\includegraphics[width=\linewidth]{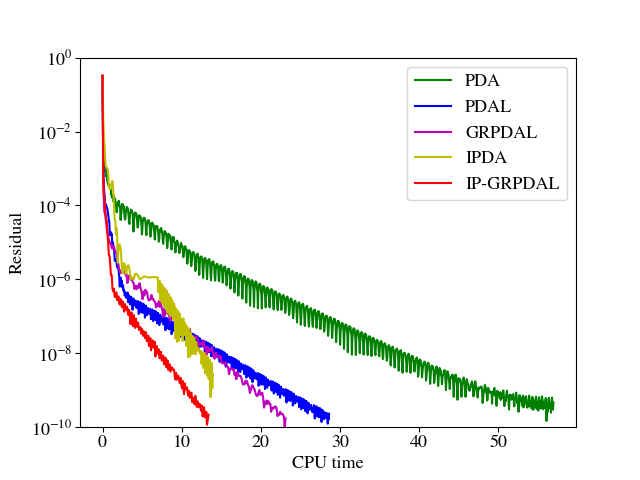}
\caption*{(a)$(n,p,s)=(100,100,10)$}
\end{minipage}
\hspace{-8pt}
\begin{minipage}[t]{0.33\textwidth}
\centering
\includegraphics[width=\linewidth]{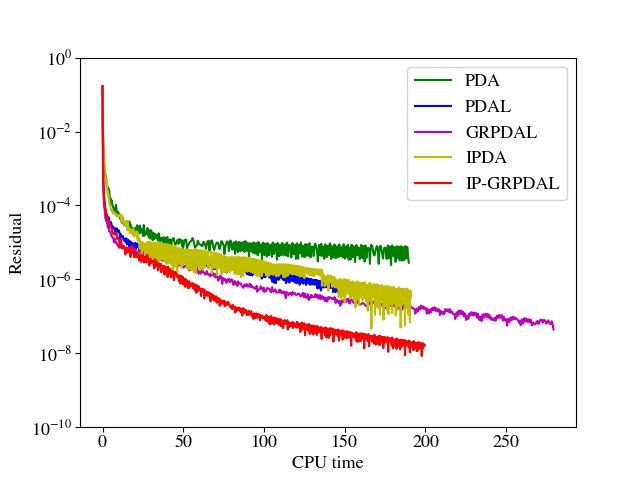}
\caption*{(b)$(n,p,s)=(500,800,50)$}
\end{minipage}
\hspace{-8pt}
\begin{minipage}[t]{0.33\textwidth}
\centering
\includegraphics[width=\linewidth]{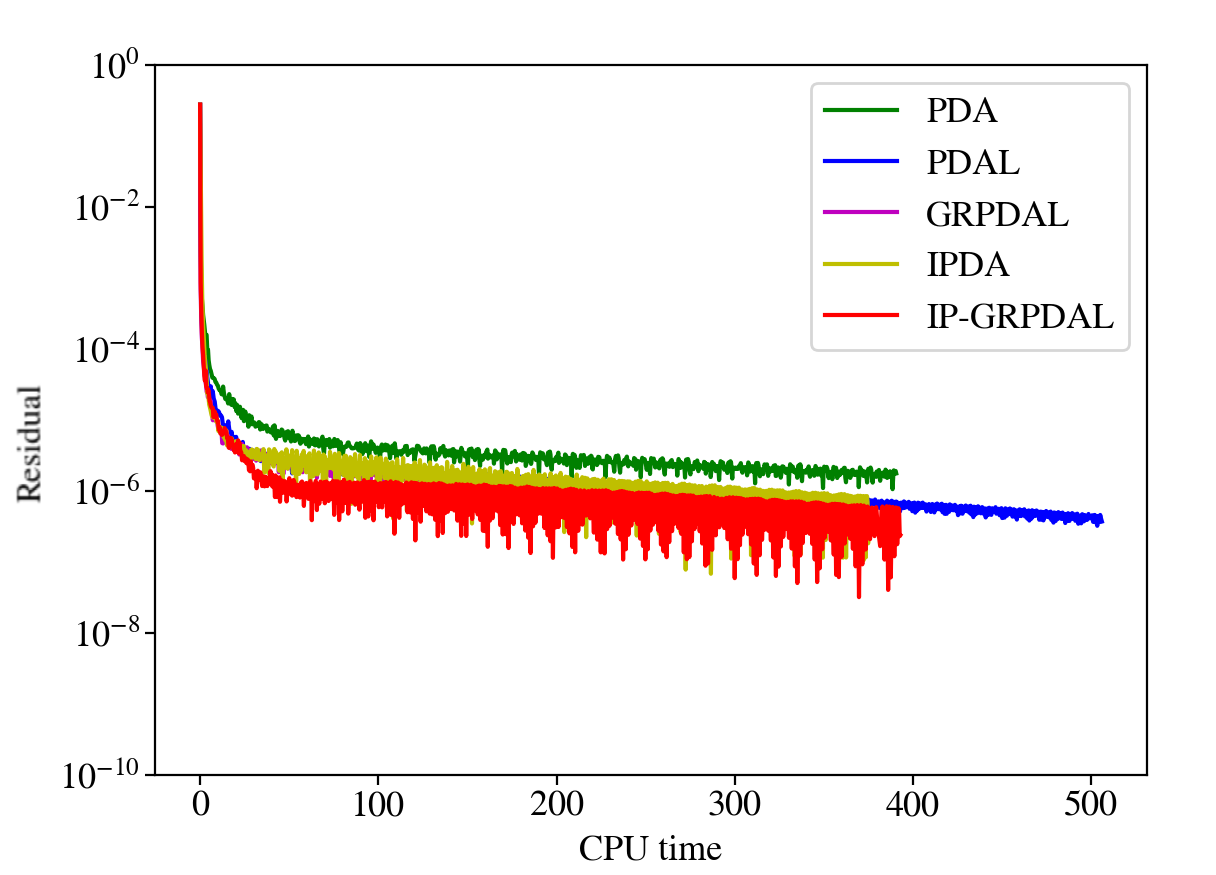}
\caption*{(c)$(n,p,s)=(1000,2000,100)$}
\end{minipage}
\captionsetup{justification=raggedright,singlelinecheck=false}
\caption{Evolution of function value residuals with respect to CPU time}
\label{fig:figure-n}
\end{figure}

From the results presented in Table \ref{tab:table-np}, we can observe that our IP-GRPDAL method takes less CPU time and number of iterations compared with the other ones. In Figure \ref{fig:figure-n}, the ordinate denotes the function value residuals $\varPhi(x^k)-\varPhi^*$ while the abscissa denotes the CPU time, from which it can be seen that the IP-GRPDAL method is much faster than the other ones.
\subsection{TV-$L_1$ image deblurring problem}
In this subsection, we study the numerical solution of the TV-$L_1$ model (4.1) in \cite{FH2022} for image deblurring
 \begin{equation} \label{4.3}
\mathop{ \textrm{min}}\limits_{x \in X }  F(x) = \|Kx-f\|_1+ \nu\|Dx\|_1,
\end{equation}
where $f \in Y$ is a given (noisy) image, $K :X \to Y$ is a known linear (blurring) operator, $D :X \to Y$ denotes the gradient operator and $\nu$ is a regularization parameter. We introduce the variables $\kappa_1,\kappa_2>0$ which satisfy $\kappa_1 + \kappa_2 = \nu$. Then, (\ref{4.3}) can be written as
$$ \mathop{ \textrm{min}}\limits_{x \in X}  \|Kx-f\|_1 + \kappa_1\|Dx\|_1+ \kappa_2\|Dx\|_1.$$
Further, the above formula can be rewritten as
$$   \mathop{ \textrm{min}}\limits_{x \in X} \mathop{ \textrm{max}}\limits_{y \in Y} L(x,y) :=   \kappa_1\|Dx\|_1 + \langle Ax,y \rangle - \Upsilon_{C_1}(u) -  \Upsilon _{C_{\kappa_2}}(v) - \langle f,u \rangle,$$
where $C_{\kappa} =\{y \in Y | \|y\|_{\infty} \leqslant \kappa\}$, $ y =\left(
 \begin{array}{c}
u \\
v
\end{array}
\right) $, and $ A =\left(
 \begin{array}{c}
K \\
\kappa_2 D
\end{array}
\right). $
We adopt the following inequality as the stopping criterion of inner loop:
$$ \varPsi(y^{k+1}, v^{k+1}) \leqslant \varepsilon_{k+1},$$
where $\varepsilon_{k+1} = O(1/(k+1)^\alpha)$ and $ \varPsi$ is defined as (4.11) in \cite{FH2022}.

\begin{figure}
\centering
\begin{minipage}[t]{0.35\textwidth}
\centering
\includegraphics[width=\linewidth]{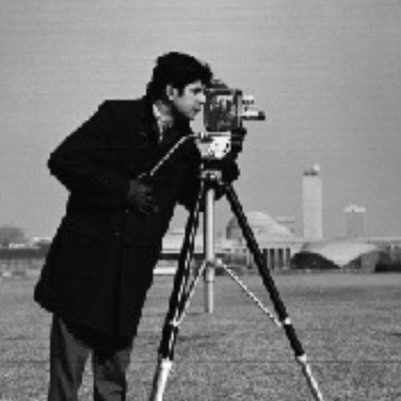}
\caption*{(a) Original Image}
\end{minipage}
\hspace{17pt}
\begin{minipage}[t]{0.35\textwidth}
\centering
\includegraphics[width=\linewidth]{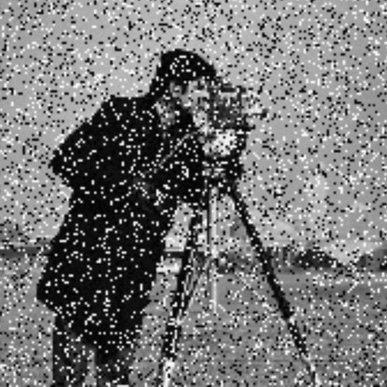}
\caption*{(b) Noise Image}
\end{minipage}
\captionsetup{justification=raggedright,singlelinecheck=false}
\caption{Cameraman.png (256 $\times$ 256)}
\label{fig:figure-c}
\end{figure}

 In the following we will report the numerical experiment results.
In this test, average blur with hsize=9 was applied to the original image cameraman.png(256 ¡Á 256) by fspecial(average, 9), and $20\%$ salt-pepper noise was added in(see Figure \ref{fig:figure-c}). At the same time, we adopt the following stopping rule:
 $$  \frac{F(x^k)-F(x^*)}{F(x^*)} < 10^{-5}, $$
 where ${x}^{*}$ is a solution of the TV-$L_1$ model (\ref{4.3}).

We fixed the number of iterations as $100$ and the penalty coefficient $ \nu = 0.1$. When the five algorithms are implemented, their respective parameters are choosen as follows:\\
\smallcircle PDA: $\tau=\sigma=0.99$;\\
\smallcircle IPDA: $\tau_0=\sigma_0=1, s_1 = 2, r_1 = \frac{0.99}{s_1}, s_2 = 1, r_2 = \frac{0.99}{s_2}, \alpha =2$;\\
\smallcircle PDAL: $\tau_0=0.1, \beta=1,\mu=0.1, \eta=0.99$;\\
\smallcircle GRPDAL: $\tau_0=0.1, \phi=1.618, \beta=1,\mu=0.1, \eta=0.99$;\\
\smallcircle IP-GRPDAL: $\tau_0=0.1, \phi=1.618, \beta=1, s_1 = 2, r_1 = \frac{0.99}{s_1}, s_2 = 1, r_2 = \frac{0.99}{s_2}, \alpha =2, \mu=0.1, \eta=0.99$.\\

\begin{figure}[!ht]
\centering
\begin{minipage}[t]{0.3\textwidth}
\centering
\includegraphics[width=\linewidth]{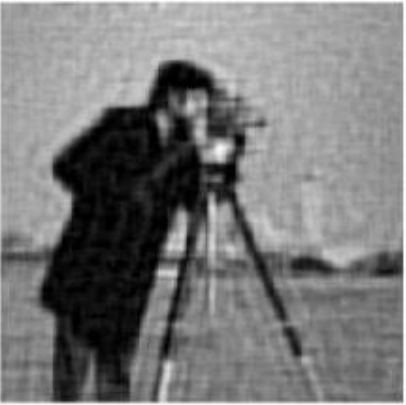}
\caption*{(a) PDA}
\end{minipage}
\hspace{10pt}
\begin{minipage}[t]{0.3\textwidth}
\centering
\includegraphics[width=\linewidth]{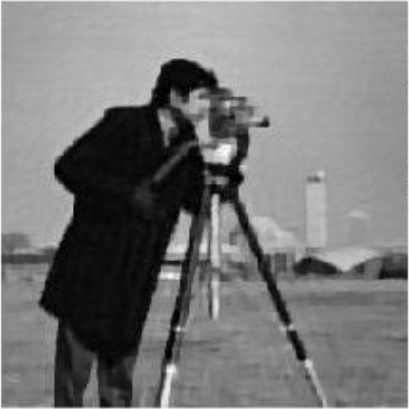}
\caption*{(b) IPDA}
\end{minipage}
\hspace{10pt}
\centering
\begin{minipage}[t]{0.3\textwidth}
\centering
\includegraphics[width=\linewidth]{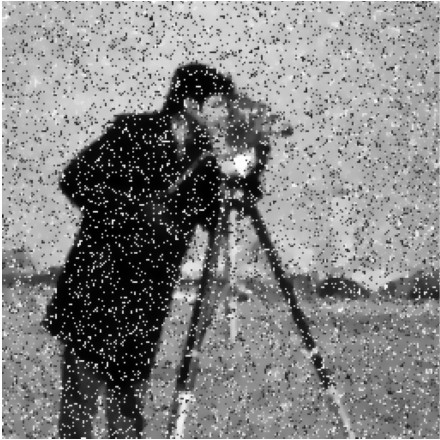}
\caption*{(c) PDAL}
\end{minipage}
\begin{minipage}[t]{0.3\textwidth}
\centering
\includegraphics[width=\linewidth]{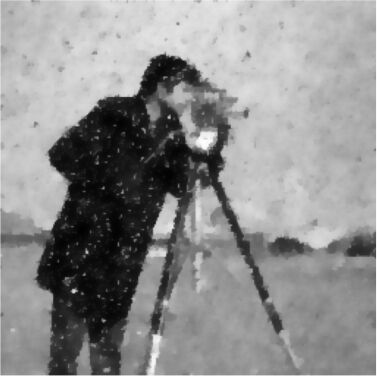}
\caption*{(d) GRPDAL}
\end{minipage}
\hspace{10pt}
\begin{minipage}[t]{0.3\textwidth}
\centering
\includegraphics[width=\linewidth]{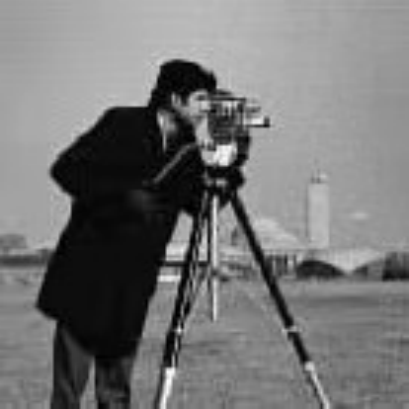}
\caption*{(e) IP-GRPDAL}
\end{minipage}
\captionsetup{justification=raggedright,singlelinecheck=false}
\caption{Restored images}
\label{fig:figure-r}
\end{figure}

The restored images by  the above Algorithms are displayed in Figure \ref{fig:figure-r}. Obviously, our  IP-GRPDAL method gets better restoration quality compared with PDA, PDAL, IPDA and GRPDAL methods. In our experiment, we find that, if we increase the number of iterations to $1000$ or more, all algorithms can restore the image with almost the same quality, but our algorithm needs fewer iterations than the other ones.

\section{Conclusions}\label{sec:clu}
In this paper, we propose an inexact golden ratio primal-dual algorithm with linesearch for the saddle point problem by applying inexact extended proximal terms with matrix norm introduced in \cite{FH2022}. Under the assumption that $\tau\sigma\|A\|_T^2 < \phi$, we show the convergence of our Algorithm \ref{A1}, provided that both controlled error sequences $\{\delta_{k+1}\}$ and $\{\varepsilon_{k+1}\}$ were required to be summable. The $O(1/N)$ convergence rate in the ergodic sense is also established in the general convex case. When either one of the component functions is strongly convex, accelerated version of Algorithm \ref{A1} is proposed, which achieves $O(1/N^2)$  ergodic convergence rate. Furthermore, the linear convergence results are established when  both component functions are strongly convex. We also apply our method to solve sparse recovery and TV-$L_1$ image deblurring problems and verify their efficiency numerically. It will be a interesting open problem to  establish nonergodic  convergence rate of  IP-GRPDAL method with linesearch.\\

\end{document}